%% file: main.tex
\documentclass{article}

\usepackage[Symbol]{upgreek}

\usepackage{amsmath,amssymb,amsfonts,mathrsfs,mathtools,dsfont,mathrsfs,bm,bbm} 
\DeclareMathAlphabet{\mathpzc}{OT1}{pzc}{m}{it}
\DeclareSymbolFontAlphabet{\amsmathbb}{AMSb}%
\usepackage{csvsimple,booktabs}
\usepackage{float,enumitem}
\usepackage[caption=false]{subfig}
\usepackage[dvipsnames]{xcolor}
\usepackage{graphicx,epsfig}
\usepackage{balance}
\usepackage[hyphens]{url}
\usepackage[breaklinks, colorlinks, linkcolor=MidnightBlue, anchorcolor=MidnightBlue, citecolor=MidnightBlue, urlcolor=MidnightBlue]{hyperref}
\usepackage{cite}

\usepackage[ruled,linesnumbered]{algorithm2e}

\usepackage{psfrag}
\usepackage{color}

\usepackage{bookmark}

\input{header}

\renewcommand{\ii}{{\bm{i}}}
\renewcommand{\bone}{{\mathds{1}}}

\newcommand{\bose}[1]{{\textcolor{purple}{Bose says: {#1}}}}

\definecolor{teal}{rgb}{0.0, 0.5, 0.5}

\newcommand{\vvec}[1]{\textrm{vec}\left(#1\right)}

\usepackage{tabularx}
\usepackage{makecell}
\newcolumntype{L}{X}
\newcolumntype{C}{>{\centering \arraybackslash}X}
\newcolumntype{R}{>{\raggedright \arraybackslash}X}
\usepackage{booktabs}
\newtheorem{theorem}{Theorem}

\begin{document}

\title{\LARGE {\bf{Distributed Dual Subgradient Methods with Averaging and Applications to Grid Optimization}}}

\author{Haitian Liu \quad Subhonmesh Bose \quad 
Hoa Dinh Nguyen \\ 
Ye Guo \quad Thinh T. Doan \quad Carolyn L. Beck
\thanks{S. Bose and C.L. Beck are with the University of Illinois at Urbana-Champaign, Urbana, IL 61801, USA. H.D. Nguyen is with the International Institute for Carbon-Neutral Energy Research (WPI-I$^2$CNER) and Institute of Mathematics for Industry (IMI), Kyushu University, 744 Motooka, Nishi-ku, Fukuoka 819-0395, Japan. H. Liu and Y. Guo are with the Tsinghua-Berkeley Shenzhen Institute, Shenzhen, Guangdong 518055, China. T.T. Doan is with Virginia Tech,  Blacksburg, VA 24060 USA. E-mails: \url{boses@illinois.edu}, \url{hoa.nd@i2cner.kyushu-u.ac.jp}, \url{liuht19@mails.tsinghua.edu.cn }, 
\url{guo-ye@sz.tsinghua.edu.cn}, \url{thinhdoan@vt.edu}, \url{beck3@illinois.edu}. This project was partially supported by grants from the Power Systems Engineering Research Center (PSERC), JSPS Kakenhi Grant Number JP23K03906, National Science Foundation of China under Grant 51977115.}}	

\maketitle

\begin{abstract}
We study finite-time performance of a recently proposed distributed dual subgradient (DDSG) method for convex constrained multi-agent optimization problems. The algorithm enjoys performance guarantees on the last primal iterate, as opposed to those derived for ergodic means for vanilla DDSG algorithms. Our work improves the recently published convergence rate of $\Ocal(\log T/\sqrt{T})$ with decaying step-sizes to $\Ocal(1/\sqrt{T})$ with constant step-size on a metric that combines suboptimality and constraint violation. We then numerically evaluate the algorithm on three grid optimization problems. Namely, these are tie-line scheduling in multi-area power systems, coordination of distributed energy resources in radial distribution networks, and joint dispatch of transmission and distribution assets. The DDSG algorithm applies to each problem with various relaxations and linearizations of the power flow equations. The numerical experiments illustrate various properties of the DDSG algorithm--comparison with vanilla DDSG, impact of the number of agents, and why Nesterov-style acceleration can fail in DDSG settings.
\end{abstract}

\input{introduction}

\input{dual_subgradient}

%
\input{multi_area}

\input{DER_coordination}
\input{TD_coordination}
{
\section{Concluding Remarks}
\label{sec:Conclusion}

We have studied a constant step-size distributed dual subgradient (DDSG) method with averaging that provides order-optimal $1/\sqrt{T}$ convergence rate after $T$ iterations for  multi-agent convex constrained optimization problems of the form $\Pcal$ whose objective function need not be strongly convex. The convergence guarantee  improves the rate of the decaying step-size, fully decentralized counterpart studied in \cite{liang2019distributed}. We formulated a variety of problems that arise in the operation of an electric power system as examples of $\Pcal$ and applied the DDSG method to solve them. The case studies with DDSG on grid optimization problems are presented not to empirically challenge other methods in the literature, but as means to illustrate the broad applicability of optimization problems such as $\Pcal$ to facilitate unified algorithm development. These examples also underscore various properties of the DDSG algorithm. Specifically, the first example illustrates the difference in performance of the last iterate between the vanilla DDSG method and ours. The second example reveals the effect of the number of agents on convergence speed. The third example empirically studies the possibility of Nesterov-style acceleration in DDSG methods and the role that smoothness in the dual function plays in such analysis.


Several research directions are of interest to us. The first among them is an algorithm with last-iterate guarantees over time-varying communication networks. The second important direction is to better disentangle the error $\Vcal_T$ to provide guarantees on suboptimality and constraint violation separately. Third, we want to explore ways to exploit a similar algorithmic architecture for stochastic and risk-sensitive convex optimization problems.

}


\appendix
\section*{Appendix}
\input{supplementary_materials}

\bibliographystyle{plain}
\bibliography{refs}


\end{document}

%% file: header.tex
\newcommand{\beq}{\begin{equation}}
\newcommand{\eeq}{\end{equation}}
\newcommand{\beqa}{\begin{eqnarray}}
\newcommand{\eeqa}{\end{eqnarray}}
\newcommand{\beqan}{\begin{eqnarray*}}
\newcommand{\eeqan}{\end{eqnarray*}}

\newcommand{\trace}{\mbox{\rm Tr}}
\newcommand{\vnorm}[1]{\left\|#1\right\|}

\newcommand{\ii}{{\bm{\mathpzc{i}}}}

\newcommand{\mean}[1]{m_{_{#1}}} 

\newcommand{\diag}{\mathop{\mathrm{diag}}}

\newcommand\T{{\mathpalette\raiseT\intercal}}
\newcommand\raiseT[2]{\raisebox{0.25ex}{$#1#2$}
}

\newcommand{\Cset}{\mathds{C}}

\newcommand{\Eset}{\mathds{E}}

\newcommand{\Rset}{\mathds{R}}
\newcommand{\Sset}{\mathds{S}}

\newcommand{\Xset}{\mathds{X}}

\newcommand{\Zset}{\mathds{Z}}

\newcommand{\Acal}{{\cal A}}

\newcommand{\Dcal}{{\cal D}}

\newcommand{\Lcal}{{\cal L}}

\newcommand{\Ncal}{{\cal N}}
\newcommand{\Ocal}{{\cal O}}
\newcommand{\Pcal}{{\cal P}}

\newcommand{\Vcal}{{\cal V}}

\newcommand{\Gfk}{{\mathfrak{G}}}

\newcommand{\Hsf}{\sf{H}}

\newcommand{\argmin}{\mathop{\rm argmin}}

\newcommand{\bone}{\mathbf{1}}

\renewcommand{\v}[1]{{\bm{#1}}}

\newcommand{\ol}[1]{\ensuremath{\overline{{#1}}}}
\newcommand{\ul}[1]{\ensuremath{\underline{{#1}}}}

\newcommand{\proj}{{\uppi}}
\newcommand{\dist}{\textrm{dist}}
\renewcommand{\mean}[1]{\ol{#1}}

\newcounter{l1}
\newcounter{l2}
\newcounter{l3}
\newcommand{\bdotlist}{\begin{list}{$\bullet$}{}}
\newcommand{\bboxlist}{\begin{list}{$\Box$}{}}
\newcommand{\bbboxlist}{\begin{list}{\raisebox{.005in}{{\tiny
$\blacksquare$ \ \ }}}{}}
\newcommand{\bdashlist}{\begin{list}{$-$}{} }
\newcommand{\blist}{\begin{list}{}{} }
\newcommand{\barablist}{\begin{list}{\arabic{l1}}{\usecounter{l1}}}
\newcommand{\balphlist}{\begin{list}{(\alph{l2})}{\usecounter{l2}}}
\newcommand{\bAlphlist}{\begin{list}{\Alph{l2}.}{\usecounter{l2}}}
\newcommand{\bdiamlist}{\begin{list}{$\diamond$}{}}
\newcommand{\bromalist}{\begin{list}{(\roman{l3})}{\usecounter{l3}}}

%% file: introduction.tex

\section{Introduction}
\label{sec:intro}

{Distributed optimization algorithms offer mechanisms to optimize in multi-agent environments, where one cannot aggregate all problem data in a central location. Agents in this paradigm iteratively perform local computational steps and communicate relevant variables over a network. While a variety of distributed solution architectures have been developed and analyzed in the literature for various settings, we focus on distributed dual subgradient (DDSG) methods with averaging that can be used to solve convex constrained multi-agent optimization problems of the form 
\begin{subequations}
\begin{alignat}{2}
\Pcal: \ 
& \underset{\v{x}_1, \ldots, \v{x}_N}{\text{minimize}} && \quad  \sum_{j=1}^N f_j(\v{x}_j), 
\label{eq:P.o}\\ 
& \text{subject to} 
&&
\quad \sum_{j=1}^N \v{g}_j^\textrm{E}(\v{x}_j) = 0,
\label{eq:P.E}
\\
&&&
\quad \sum_{j=1}^N \v{g}_j^\textrm{I}(\v{x}_j) \leq 0,
\label{eq:P.I}
\\
&&& \quad  \v{x}_j \in \Xset_j \subseteq \Rset^{n_j}, \ j= 1,\ldots, N.
\label{eq:P.X}
\end{alignat}
\label{eq:P}
\end{subequations}
The $N$ agents communicate only across edges of an undirected graph $\Gfk(N, \Eset)$. Assume that functions $f_j$, $\v{g}_j^\textrm{I}$ are convex and  $\v{g}_j^\textrm{E}$ is affine over the compact convex set $\Xset_j$ for each $j=1,\ldots,N$.

At its core, DDSG methods rely on dual decomposition that starts by separating the Lagrangian into agent-wise Lagrangian functions that each agent optimizes, given a dual iterate (see classical texts such as \cite{BertsekasTsitsiklis1989,polyak1987introduction}). This agent-wise optimization of the primal variables, given a dual iterate, can be shown to provide a subgradient of the dual function at the dual iterate. Thus, a subgradient ascent on the dual function can be achieved to solve the dual problem. Such an update rule requires a central coordinator to manage the dual iterates, which is undesirable in many distributed contexts.  DDSG methods circumvent this need by maintaining \emph{local copies} of such multipliers and running a consensus-based distributed dual ascent on these local multiplier copies. Approximate primal solutions can be recovered from these dual solutions as in \cite{simonetto2016primal}, building on techniques in \cite{Nedic2009,gustavsson2015primal,ma2007recovery,larsson1999ergodic}, among others; asymptotic guarantees on recovered primal sequences are also known.

In this paper, we study a distributed dual subgradient method, analyzed recently in \cite{liang2019distributed}, that provides a fully distributed variant of the algorithm, also proposed recently in \cite{nesterov2018dual}. As opposed to vanilla DDSG algorithms, convergence guarantees of these algorithms are obtained for the \emph{last primal iterate} as opposed to that for \emph{ergodic means} of primal sequences. As Nesterov and Shikhman argue in \cite{nesterov2018dual}, iterates in vanilla DDSG algorithms can oscillate wildly during algorithm execution. If the intermediate primal iterates of the algorithm are implemented in practice, vanilla DDSG can negatively impact stability in multi-agent environments. Besides, primal recovery is no longer required for asymptotic constraint satisfaction in this paradigm.

For the algorithms in \cite{liang2019distributed,nesterov2018dual}, 
asymptotic guarantees for the convergence of the last primal-dual iterate to the set of primal-dual optimizers have been established with decaying step-sizes. Also, the authors of \cite{liang2019distributed} characterize a $\Ocal(\log T/\sqrt{T})$-bound on a metric that combines sub-optimality and constraint violation with decaying step-size (similar to \cite{duchi2011dual}), while they argue a lower bound of $\Ocal(1/\sqrt{T})$. Their analysis shows an $\Ocal(1/\sqrt{T})$-bound, when a coordinator manages the dual updates centrally. Thus, they view the $\log T$ factor as a price for decentralization. Our first contribution of this paper is to close this gap in Section \ref{sec:algorithm} with the proof in Section \ref{sec:proof}, i.e., we sharpen the convergence rate to $\Ocal(1/\sqrt{T})$ for this fully distributed algorithm with a \emph{constant step-size}. The dual function for $\Pcal$ can be nonsmooth, and as a result, a first-order algorithm that climbs the dual function cannot converge at a better rate (see Appendix \ref{sec:dual_func}), making the algorithm order-optimal.

Through the rest of the paper, we study various properties of our DDSG algorithm on $\Pcal$. Specifically, we consider three different optimization problems that arise in operations of the electric power system  and cast them as examples of $\Pcal$. 
For each problem, we adopt different power flow models, different notions of agents and different definitions of the network over which these agents interact. This exercise stands as our second contribution, which demonstrates the modeling power of $\Pcal$ and the broad applicability of the algorithm with convergence guarantees. The grid optimization problems are nonconvex due to the nature of Kirchhoff's laws \cite{cain2012history}. In this paper, we \emph{convexify} each problem by considering approximations and convex relaxations for power flow models that are suited to that application. Each problem and the application of DDSG illustrates specific properties of the algorithm that are delineated below.
    
    $\bullet$ \emph{Multi-area optimal power flow  ($\Pcal_1$)} seeks to dispatch resources over an interconnected transmission network, parts of which are controlled by different system operators. The distributed algorithmic architecture we study in Section \ref{sec:multi_area}, bypasses the need for the system operators to share all relevant data from within their footprint with another system operator and yet seek to solve a joint optimal power flow problem through a distributed solution architecture, e.g., see \cite{Guo2018}. For $\Pcal_1$, we consider a linear power flow model for transmission networks motivated by the use of such models in wholesale market environments that facilitate inter-area coordination, e.g., in \cite{stott2009dc}. Through the numerical example on $\Pcal_1$, we demonstrate how the averaging scheme stabilizes the last iterate, compared to vanilla DDSG methods, and circumvents the need for primal recovery. We also show that finite-time constraint violation of the algorithm is empirically better than its theoretical bound suggests.

    $\bullet$ \emph{Coordination of distributed energy resources (DERs) in distribution grids ($\Pcal_2$)} is designed to optimize real and reactive power outputs from DERs at the grid-edge to minimize cost (dis-utility) of such a dispatch and possibly a network-wide objective such as frequency regulation. A distributed algorithm allows quick updates of optimization variables without the need to communicate with a central coordinator across the distribution grid. For $\Pcal_2$ in Section \ref{sec:DER_coord}, we consider a second-order cone programming (SOCP) based relaxation of the power flow equations in the distribution grids. The thoroughly-studied SOCP-based relaxation of power flow equations for distribution grids are often tight in practice (see \cite{Farivar2013,gan2015exact}). Through examples of $\Pcal_2$, we show that our algorithm can track changing grid conditions for small networks. Tracking performance degrades with network size, where aggregation of nodes as agents becomes vital for performance. Said differently, if speed is paramount, one must carefully control the degree of decentralization.
    
     $\bullet$ \emph{Transmission and distribution (T\&D) grid coordination ($\Pcal_3$)} seeks to dispatch assets across the transmission and distribution grids without the need to collect all information from the grid-edge and the bulk power systems at one location. The distributed solution architecture in Section \ref{sec:TD_coordination} alleviates the transmission system operator's lack of visibility into utility-managed distribution networks. For $\Pcal_3$, we consider two different power flow models for the transmission and the distribution grids. For transmission, we choose a semidefinite programming (SDP) based relaxation of power flow equations, given its popularity to approach the AC optimal power flow problem \cite{lavaei2011zero,zhang2012geometry,bose2015quadratically}. For the distribution grids, we consider the linear distribution power flow model from \cite{Baran1989c}.
    Using an example of $\Pcal_3$, we compare the empirical performance of the DDSG method with an adopted distributed acceleration scheme proposed in \cite{qu2019accelerated}, where acceleration fails to increase convergence speed with linear costs, but performs well with quadratic costs. We argue  that nonsmoothness of the dual function with linear costs is behind this phenomenon, and elaborate on its role in Appendix \ref{sec:dual_func}.
}   

We recognize that a long literature has emerged on each of the grid optimization problems; we only cite a few, owing to space constraints. Our goal in presenting the power system examples is not an attempt to empirically compare our DDSG algorithm with others in the literature proposed for these problems. Rather, our formulations of these problems as instances of $\Pcal$ reveal the generality of $\Pcal$ and make it possible for subsequent unified algorithm development for $\Pcal$. In addition, each example highlights an important aspect of the DDSG algorithm.

%% file: dual_subgradient.tex
\section{The dual subgradient method with averaging}
\label{sec:algorithm}
We present the DDSG algorithm to solve $\Pcal$ in \eqref{eq:P}.
The finite-time performance guarantee for $\Pcal$ is derived in Section \ref{sec:proof}. Then in Sections \ref{sec:multi_area}, \ref{sec:DER_coord} and \ref{sec:TD_coordination}, we cast $\Pcal_1$, $\Pcal_2$ and $\Pcal_3$ as examples of $\Pcal$. 

To lighten notation, let $\v{g}_j$ collect both $\v{g}_j^\textrm{E}$ and $\v{g}_j^\textrm{I}$ with the understanding that the first $M^\textrm{E}$ constraints encode equalities and the last $M^\textrm{I}$ are inequalities.
The algorithm relies on Lagrangian duality theory associated with $\Pcal$. We begin by defining the Lagrangian function
\begin{align}
\begin{aligned}
    \Lcal(\v{x}, \v{z}) &:=  
    \sum_{j=1}^N \left[ f_j(\v{x}^j) +  \v{z}^{\T} \v{g}_{j} (\v{x}_j) \right], 
\end{aligned}
\label{eq:defL}
\end{align}
\begin{gather}
\begin{gathered}
    \v{x}^\T = (\v{x}_1^\T, \ldots, \v{x}_N^T) \in \Xset:= \Xset_1 \times \ldots \times \Xset_N,\\
    \v{z} \in \Zset := \Rset^{M^{\textrm{E}}} \times \Rset^{M^\textrm{I}}_+.
\end{gathered}
\end{gather}
Then, $\Pcal$ can be cast as a min-max problem with optimal value $\Pcal^\star$ as
\begin{align}
    {\Pcal}^\star = \min_{\v{x} \in \Xset} \max_{\v{z} \in \Zset} \  \Lcal(\v{x}, \v{z}).
    \label{eq:P.primal}
\end{align}
Let $\Xset^\star$ denote the set of optimizers of $\Pcal$. Associated with $\Pcal$ is its dual problem
\begin{align}
    \Pcal_D^\star = \max_{\v{z} \in \Zset} \min_{\v{x} \in \Xset} \  \Lcal(\v{x}, \v{z}).
    \label{eq:P.dual}
\end{align}
Let $\Zset^\star$ denote the set of optimizers of the dual problem. 
Weak duality implies that $\Pcal^\star \geq \Pcal_D^\star$. We say strong duality holds if the inequality is met with an equality. And, $\v{x}^\star \in \Xset$, $\v{z}^\star \in \Zset$ is a saddle point of $\Lcal$, if 
\begin{align}
    \Lcal(\v{x}^\star, \v{z}) \leq \Lcal(\v{x}^\star, \v{z}^\star) \leq \Lcal(\v{x}, \v{z}^\star),
    \label{eq:saddle}
\end{align}
for all $\v{x} \in \Xset, \v{z} \in \Zset$. The well-known saddle point theorem (see \cite[Theorem 2.156]{bonnans2013perturbation}) states that the  primal-dual optimizers $\Xset^\star \times \Zset^\star$ coincide with the saddle points of $\Lcal$. 
%
We assume throughout that the set of saddle points of $\Pcal$ is nonempty and bounded. As a result, strong duality holds for $\Pcal$, i.e., $\Pcal^\star = \Pcal_D^\star$, and the set of primal-dual optimizers is nonempty. Saddle-points exist under standard constraint qualifications such as Slater's condition, e.g., see \cite[Theorem 2.165]{bonnans2013perturbation}.

Dual decomposition techniques for distributed optimization rely on the observation that the dual function separates into agent-wise optimization problems, given a multiplier $\v{z}$ as
\begin{align} 
 \min_{\v{x} \in \Xset} \ \Lcal(\v{x}, \v{z}) = \sum_{j=1}^N \underbrace{\min_{\v{x}_j \in \Xset_j}\Lcal_j (\v{x}_j, \v{z})}_{:= \Dcal_j(\v{z})},
 \label{eq:dual.decomp}
\end{align}
where $\Lcal_j(\v{x}_j, \v{z}) :=  f_j(\v{x}^j)
    + \v{z}^{\T} \v{g}_j(\v{x}_j)$.
If the agents can perform these agent-wise minimizations, then  a distributed projected subgradient ascent algorithm can solve the dual problem (e.g., see \cite{boyd2007notes}). Per Danskin's theorem, a subgradient $\nabla_\v{z}\Dcal_j(\v{z})$ can be obtained from the agent-wise minimization of $\Lcal_j$, given that the sub-differential set of the concave function $\Dcal_j$ at $\v{z}$ is
\begin{align}
    \partial_z \Dcal_j(\v{z}) := \textrm{conv}\{ \partial_z \Lcal_j(\v{x}_j, \v{z}) \ \mid \ \v{x}_j \in \Xset_j^\star(\v{z})  \}.
    \label{eq:diff.Dj}
\end{align}
Here, ``$\textrm{conv}$'' computes the convex hull of its argument and $\Xset_j^\star(\v{z})$ is the set of minimizers of $\Lcal_j(\cdot, \v{z})$ over $\Xset_j$. The minimization problem is well-defined, given that $\Xset_j$'s are  compact. 
Running such an algorithm, however, requires a \emph{central coordinator} to compute the $\v{z}$-update and broadcast the results to all agents. Albeit simpler than aggregating all problem data at a single location, the need for said coordination is a downside of classical dual decomposition. 

To avoid coordination for the dual update, one can alternately create \emph{local copies} of $\v{z}$'s among all agents and enforce equality among these local estimates in the dual problem as
\begin{align}
\begin{aligned}
\max \ \sum_{j=1}^N \Dcal_j(\v{z}_j),
 \text{subject to }  \v{z}_{j} = \v{z}_{k}, \ j,k=1,\ldots,N,
 \end{aligned}
 \label{eq:P.dual.local}
\end{align} 
where $\v{z}_j$ is the local copy of $\v{z}$ with agent $j$. One can run a projected distributed subgradient ascent as in \cite{simonetto2016primal} to solve  \eqref{eq:P.dual.local}. The primal iterates obtained from agent-wise minimization of $\Lcal_j$ evaluated at the dual iterates may fail to collectively satisfy the constraints of $\Pcal$. Primal averaging schemes have been studied in \cite{simonetto2016primal}; limit points of such recovered primal solutions are known to satisfy the constraints. One can judiciously maintain local copies only among a subset of the agents to relieve communication burden (see \cite{kao2019convergence}).

Recently, a dual subgradient algorithm was proposed in \cite{nesterov2018dual} that leveraged an {estimation sequence} technique to provide guarantees on sub-optimality and infeasibility on the \emph{last iterate}. This algorithm does not treat ergodic means simply as outputs from a dual subgradient calculation, but rather uses these means as primal-dual iterates to run the algorithm. We focus on the fully distributed variant of the algorithm that is proposed and analyzed in \cite{liang2019distributed}. To present the algorithm, let $\v{W} \in \Rset^{N\times N}$ be a doubly stochastic, irreducible and aperiodic weighting matrix that follows the sparsity pattern of $\Gfk$, i.e.,
\begin{align}
    W_{j,k} \neq 0 \iff (j, k) \in \Eset.
\end{align}
Then, the distributed projected dual subgradient with averaging is given by Algorithm \ref{alg:ddsa}, where $\v{x}_j/\v{X}_j$ are primal sequences and $\v{z}_j/\v{Z}_j$ are dual sequences. The updates comprise minimization of the local dual function in step \ref{alg:X}, averaging of these primal minimizers in step \ref{alg:x}, a consensus followed by local subgradient-based dual update in step \ref{alg:Z} and an ergodic mean computation for the projected dual variable in step \ref{alg:z} with step-size $\eta$. Here, $\proj_\Zset$ projects the arguments on $\Zset$.


\begin{algorithm}
	\caption{Distributed dual subgradient with averaging to solve $\Pcal$.}
	\label{alg:ddsa}

    Choose $\v{z}_j(1) = 0$, $\v{Z}_j(0) = 0$, $\v{x}_j(0) \in \Xset_j$ and $\eta = \eta_0/\sqrt{T}$.

	\For{$t=1, \ldots, T$}{
	
	    $\v{X}_j(t) \gets \argmin_{\v{x}_j \in \Xset_j} \ \Lcal_j(\v{x}_j, \v{z}_j(t))$. \label{alg:X}

	    $\v{x}_j(t) \gets \frac{t-1}{t} \v{x}_j(t-1) + \frac{1}{t} \v{X}_j(t)$. \label{alg:x}

        $\v{Z}_j(t) \gets \sum_{k=1}^N W_{jk} \v{Z}_k(t-1)+ t \v{g}_j(\v{x}_j(t)) - (t-1) \v{g}_j(\v{x}_j(t-1))$.\label{alg:Z}

        $\v{z}_j({t+1})\gets  \frac{t}{t+1}\v{z}_j(t) + \frac{1}{t+1} \proj_\Zset \left[ \eta \v{Z}_j(t) \right]$. \label{alg:z}
    
    }
\end{algorithm}

To study convergence properties of this algorithm, consider the metric introduced in \cite{nesterov2018dual} and used in \cite{liang2019distributed}, given by 
\begin{align}
\begin{aligned}
    \mathscr{V}_T\left(\v{x}(T), \mean{\v{z}}(T)\right) 
    := \sum_{j=1}^N f_j(\v{x}_j(T)) - \sum_{j=1}^N \Dcal_j(\mean{\v{z}}(T)) 
     + \frac{\eta T}{2N} \vnorm{\proj_\Zset \left[ \sum_{j=1}^N \v{g}_j(\v{x}_j(T)) \right]}^2,
\end{aligned}
\label{eq:V.def}
\end{align}
where 
$\mean{\v{z}}(T) := \frac{1}{T} \sum_{t=1}^T \proj_\Zset \left[ \eta \mean{\v{Z}}(t-1) \right]$, $\mean{\v{Z}}(t) := \frac{1}{N} \sum_{j=1}^N \v{Z}_j(t)$.
The sum of the first two terms measures the gap between the primal objective at $\v{x}(T) \in \Xset$ and the dual function evaluated at $\mean{\v{z}}(T) \in \Zset$. The last summand is a measure of the constraint violation at $\v{x}(T)$. We sharpen the bound of \cite[Theorem 2]{liang2019distributed} in the next result. 
\begin{theorem}
\label{thm:main}
Iterates generated by Algorithm \ref{alg:ddsa} with $\eta = \eta_0/\sqrt{T}$ over $t=1, \ldots, T$, $\eta_0>0$ constant, satisfy
\begin{align}
\begin{aligned}
        \mathscr{V}_{T}(\v{x}(T), \mean{\v{z}}(T)) 
        &\leq \frac{1}{\sqrt{T}}\left( \frac{C_0}{1-\sigma_2(\v{W})} + C_1 \right),
        \\
        \mathscr{V}_{T}(\v{x}(T), \mean{\v{z}}(T)) 
        &\geq \Pcal^\star - \sum_{j=1}^N \Dcal_j(\mean{\v{z}}(T))  - \frac{C_2}{\sqrt{T}}.
\end{aligned}
\label{eq:result}
\end{align}
where $C$'s are positive constants that do not depend on $\Gfk$ or $T$, and $\sigma_2(\v{W})$ is the second largest singular value of $\v{W}$.
\end{theorem}
Our upper bound in this result sharpens the conclusion of \cite[Theorem 2]{liang2019distributed}, while the lower bound is identical. The result implies that the metric in \eqref{eq:V.def} indeed converges at a rate of $1/\sqrt{T}$. Our proof of the bounds largely mirrors that of \cite[Theorem 2]{liang2019distributed}, but deviates from the reliance on results from \cite{duchi2011dual} that incur the $\log T$ factor. Instead, we use an argument inspired by the proof of \cite[Theorem 2]{doan2018convergence}.

We briefly remark on the implication of Theorem \ref{thm:main} on the sub-optimality of $\v{x}(T)$ and the constraint violation, separately. Call the right hand side of the upper bound in \eqref{eq:result} as $C'/\sqrt{T}$. Then, we infer
\begin{align}
\begin{aligned}
        \sum_{j=1}^N f_j(\v{x}_j(T))  - \Pcal^\star
        \leq \sum_{j=1}^N f_j(\v{x}_j(T)) - \sum_{j=1}^N \Dcal_j(\mean{\v{z}}(T))
        \leq \frac{C'}{\sqrt{T}},
\end{aligned}
\label{eq:result.opt.gap}
\end{align}
since $\Pcal^\star$ dominates the dual function, and the second summand of $\mathscr{V}_T$ in \eqref{eq:V.def} is non-negative. 
Also, combining the two inequalities in \eqref{eq:result}, we get 
\begin{align}
\begin{aligned}
        \Pcal^\star - \sum_{j=1}^N \Dcal_j(\mean{\v{z}}(T))
        \leq \frac{C' + C_2}{\sqrt{T}}.
\end{aligned}
\label{eq:result.opt.gap.2}
\end{align}
Thus, the last primal and the dual iterate exhibit an $\Ocal(1/\sqrt{T})$ sub-optimality. The dual function can be nonsmooth at an optimum. This convergence rate is therefore order-optimal, given \cite[Theorem 3.2.1]{Nesterov04}. { Faster convergence guarantees through a Nesterov-style acceleration require stronger assumptions. See Section \ref{sec:TD_coordination} and Appendix \ref{sec:dual_func} for a discussion.}

Providing the same rate for constraint violation using Theorem \ref{thm:main} remains challenging. The difficulty stems from the fact that, unless $\v{x}(T)$ is feasible in $\Pcal$, the primal-dual gap can assume negative values. However, this gap is bounded below. Using \eqref{eq:result.opt.gap}, we obtain
\begin{align}
\begin{aligned}
\sum_{j=1}^N f_j(\v{x}_j(T)) - \sum_{j=1}^N \Dcal_j(\mean{\v{z}}(T)) 
\geq \min_{\v{x} \in \Xset} \sum_{j=1}^N f_j(\v{x}_j) - \max_{\v{x} \in \Xset} \sum_{j=1}^N f_j(\v{x}_j)
=: -D_f.
\end{aligned}
\end{align}
The constant $D_f\geq 0$ is finite, owing to the compact nature of $\Xset$. Then, \eqref{eq:result} implies
\begin{align}
\begin{aligned}
        \frac{\eta_0}{2N}\vnorm{\proj_\Zset \left[ \sum_{j=1}^N \v{g}_j(\v{x}(T)) \right]}^2
        &\leq \frac{D_f}{\sqrt{T}} + \frac{C'}{T}, 
\end{aligned}
\label{eq:result.const.viol}
\end{align}
This suggests a worst-case $\Ocal(T^{-1/4})$ decay in constraint violation--an estimate that is overly conservative as our numerical estimates will reveal. Better guarantees for vanilla DDSG methods are known, e.g., in \cite{NedicO2009_PD, Nedic2009}. A constant step-size of $\eta_0/\sqrt{T}$ yields an $\Ocal(1/\sqrt{T})$ convergence of the ergodic mean of the primal iterates.

{ With non-summable and square-summable decaying step-sizes, vanilla DDSG methods converge to a \emph{single} dual optimizer (not just to the optimal set), even in distributed settings, e.g., see \cite{gustavsson2015primal}. While asymptotic convergence to the primal-dual optimal set for Algorithm \ref{alg:ddsa} is established in \cite[Lemmas 1, 2]{liang2019distributed}, convergence of the dual iterates to a single dual optimizer has not been established. It remains unclear whether such a result is attainable beyond vanilla DDSG; the conclusion does not hold even for centralized first-order primal-dual methods; see \cite[Section 3.2]{madavan2021stochastic}.}

\input{newProof}


%% file: newProof.tex

\section{Proof of Theorem \ref{thm:main}}
\label{sec:proof}
We begin by defining additional notation.
Since $f_j$ and $\v{g}_j$'s are convex (and hence, continuous) and $\Xset_j$ is compact for each $j$, these functions admit positive constants $D_\Xset$, $D_G$, $L_g$ such that
\begin{gather}
{ 
\begin{gathered}
\vnorm{\v{x}_j - \v{x}'_j} \leq D_\Xset,    
    \quad 
    \vnorm{\v{g}_j(\v{x}_j)} \leq D_g,
    \quad
    \vnorm{\v{g}_j(\v{x}_j) - \v{g}_j(\v{x}'_j)} \leq L_g \vnorm{\v{x}_j - \v{x}'_j}
\end{gathered}
}
\end{gather}
for all $\v{x}_j, \v{x}'_j \in \Xset_j$. Also, we define 
    $D_Z := L_g D_\Xset + D_g$.

\subsection{Upper bounding $\mathscr{V}_T$}
Using this notation, we derive the upper bound on $\mathscr{V}_T$ in four steps:
\begin{enumerate}[label=(\alph*)]
    \item We bound the optimality gap as 
    \begin{align}
    {
    \begin{aligned}
        \sum_{j=1}^N \left[ f_j(\v{x}_j(T))  - \Dcal_j(\mean{\v{z}}(T)) \right]
        &\leq \frac{2D_g}{T} \sum_{j=1}^N 
     \sum_{t=1}^T   \eta \vnorm{ \v{Z}_j(t-1) -  \mean{\v{Z}}(t-1)}
    \\
    & \qquad -  \frac{\eta }{T} \sum_{j=1}^N  \sum_{t=1}^T \v{g}_j(\v{X}_j(t))^\T   \proj_\Zset \left[  \mean{\v{Z}}(t-1)  \right],
    \end{aligned}
    \label{eq:opt.gap}
    }
    \end{align}
    \item Then, we bound the constraint violation as
    \begin{align}
    {
    \begin{aligned}
    \frac{T}{2N} \vnorm{\proj_\Zset \left[ \sum_{j=1}^N \v{g}_j(\v{x}(T)) \right]}^2 
    \leq  \frac{1}{T}\sum_{t=1}^T  \sum_{j=1}^N  \v{g}_j\left(  \v{X}_j(t) \right)^\T \proj_\Zset \left[ \mean{\v{Z}}(t-1)  \right] + \frac{1}{2} N D_Z^2.
    \end{aligned}
    \label{eq:const.viol}
    }
\end{align}
    
    \item We prove that $\v{Z}_j$'s remain close to their centroid as
    \begin{align}
    {
        \sum_{j=1}^N \vnorm{ \v{Z}_j(t) -  \mean{\v{Z}}(t)}_2
        \leq  N^{3/2} D_Z \left( 1-\sigma_2(\v{W}) \right)^{-1}. 
        \label{eq:Z.centroid}
    }
    \end{align}
    
    \item Steps (a), (b), (c) are combined to prove the result.
\end{enumerate}


\noindent $\bullet$ \emph{Step (a). Bounding the duality gap:} Note that 
\begin{align}
    {
    \begin{aligned}
        \sum_{j=1}^N \left[f_j(\v{x}_j(T)) 
        - \Dcal_j(\mean{\v{z}}(T)) \right]
        &= \sum_{j=1}^N \left[ f_j(\v{x}_j(T)) - \Dcal_j(\v{z}_j(T)) \right]
        \\
        & \qquad + \sum_{j=1}^N \left[
        \Dcal_j(\v{z}_j(T)) - \Dcal_j(\mean{\v{z}}(T)) \right]
        \\
        &\leq \sum_{j=1}^N \underbrace{\left[ \frac{1}{T} \sum_{t=1}^T f_j\left(\v{X}_j(t)\right) - \Dcal_j(\v{z}_j(T)) \right]}_{:=\mathscr{A}_j}
        \\
        & \qquad + D_g \sum_{j=1}^N 
         \underbrace{\vnorm{\v{z}_j(T) - \mean{\v{z}}(T)}}_{:=\mathscr{B}_j}.
    \end{aligned}
    \label{eq:opt.gap.1}
    }
\end{align}
The last line follows from three observations: $f_j$ is convex, $\v{x}_j(T) = \frac{1}{T}\sum_{t=1}^T \v{X}_j(t)$ and $\Dcal_j$
is $D_g$-Lipschitz. In the rest of step (a), we individually bound $\mathscr{A}_j$ and $\mathscr{B}_j$.

To obtain a bound on $\mathscr{A}_j$, note that
\begin{align}
{
    t \v{z}_j(t) - (t-1) \v{z}_j(t-1)
    = \proj_\Zset\left[ \eta \v{Z}_j(t-1)\right],
}
\end{align}
which then implies
\begin{align}
    {
    \begin{aligned}
    t \Lcal_j(\v{X}_j(t), \v{z}_j(t))
    &= \Lcal_j(\v{X}_j(t), t\v{z}_j(t) - (t-1)\v{z}_j(t-1))
    \\
    & \qquad 
    + (t-1)\Lcal_j(\v{X}_j(t), \v{z}_j(t-1))
    \\
    & \geq \Lcal_j\left(\v{X}_j(t), \proj_\Zset\left[ \eta \v{Z}_j(t-1)\right]\right)
    \\
    &\qquad + (t-1)\Lcal_j(\v{X}_j(t-1), \v{z}_j(t-1)).
    \end{aligned}
    }
\end{align}
The first line follows from elementary algebra, while the second line requires the definition of $\v{Z}_j$ and the fact that $\v{X}_j(t-1)$ minimizes $\Lcal_j(\cdot, \v{z}_j(t-1))$ over $\Xset_j$. Iterating the above inequality, we obtain
\begin{align}
    {
    \begin{aligned}
    T\Dcal_j(\v{z}_j(T))
    =T \Lcal_j(\v{X}_j(T), \v{z}_j(T))
    \geq \sum_{t=1}^T \Lcal_j\left(\v{X}_j(t), \proj_\Zset\left[ \eta \v{Z}_j(t-1)\right]\right).
    \end{aligned}
    }
\end{align}
The above relation bounds $\mathscr{A}_j$ from above as
\begin{align}
    {
    \begin{aligned}
    \mathscr{A}_j 
    &\leq \frac{1}{T} \sum_{t=1}^T \left[ f_j\left(\v{X}_j(t)\right) 
    -  \Lcal_j\left(\v{X}_j(t), \proj_\Zset\left[ \eta \v{Z}_j(t-1)\right]\right) \right]
    \\
    &=-\frac{1}{T} \sum_{t=1}^T \v{g}_j(\v{X}_j(t))^\T   \proj_\Zset \left[ \eta \v{Z}_j(t-1) \right]
    \\
    &=
    - \frac{1}{T} \sum_{t=1}^T \v{g}_j(\v{X}_j(t))^\T   \left(\proj_\Zset \left[ \eta \v{Z}_j(t-1) \right] - \proj_\Zset\left[\eta \mean{\v{Z}}(t-1) \right]\right)
    \\
    & \qquad - \frac{1}{T} \sum_{t=1}^T \v{g}_j(\v{X}_j(t))^\T   \proj_\Zset \left[ \eta  \mean{\v{Z}}(t-1)  \right].
    \end{aligned}
    }
\end{align}
Cauchy-Schwarz inequality and the bounded nature of $\v{g}_j$ imply
\begin{align}
    {
    \begin{aligned}
    \mathscr{A}_j 
    &\leq
    \frac{D_g}{T} \sum_{t=1}^T  \vnorm{\proj_\Zset \left[ \eta \v{Z}_j(t-1) \right] - \proj_\Zset\left[\eta \mean{\v{Z}}(t-1) \right]}
    \\
    & \qquad - \frac{\eta}{T} \sum_{t=1}^T \v{g}_j(\v{X}_j(t))^\T   \proj_\Zset \left[  \mean{\v{Z}}(t-1)  \right].
    \end{aligned}
    \label{eq:Aj.bound}
    }
\end{align}

To bound $\mathscr{B}_j$, we use the definition of $\v{z}_j(t)$ to infer
\vspace*{-.05in}
\begin{align}
{
    \v{z}_j(T) = \frac{1}{T}\sum_{t=1}^{T} \proj_\Zset\left[\eta \v{Z}_j(t-1) \right],
}
\end{align}
which in turn implies
\begin{align}
    {
    \mathscr{B}_j 
    \leq \frac{1}{T} \sum_{t=1}^{T} \vnorm{\proj_\Zset\left[\eta \v{Z}_j(t-1) \right] - \proj_\Zset\left[\eta \mean{\v{Z}}(t-1) \right]}.
    \label{eq:Bj.bound}
    }
\end{align}
Using the bounds of \eqref{eq:Aj.bound} and \eqref{eq:Bj.bound} in \eqref{eq:opt.gap.1} and appealing to the non-expansive nature of the projection operator yields \eqref{eq:opt.gap}, completing step (a) of the proof.

\noindent $\bullet$ \emph{Step (b). Bounding the constraint violation:}
From the $\v{Z}$-update, we obtain
\vspace{-.05in}
\begin{align}
{
    \mean{\v{Z}}(t)
    =  \frac{t}{N} \sum_{j=1}^N  \v{g}_j(\v{x}_j(t)),
}
\end{align}
that proves useful in bounding the constraint violation as
\vspace{-.05in}
\begin{align}
{
\begin{aligned}
   \frac{T^2}{N^2} \vnorm{\proj_\Zset \left[ \sum_{j=1}^N \v{g}_j(\v{x}_j(T)) \right]}^2 
   &= \vnorm{\proj_\Zset \left[  \mean{\v{Z}}(T) \right]}^2
   \\
   &= \sum_{t=1}^T {\left( \vnorm{\proj_\Zset \left[ \mean{\v{Z}}(t)  \right]}^2 -  \vnorm{\proj_\Zset \left[ \mean{\v{Z}}(t-1)  \right]}^2 \right)}
   \\
   &\leq  2 \sum_{t=1}^T 
   \underbrace{\left[ \proj_\Zset \left[ \mean{\v{Z}}(t-1)  \right]^\T [\mean{\v{Z}}(t) - \mean{\v{Z}}(t-1)]
   \right]}_{:=\mathscr{E}(t)} 
   \\
   & \qquad 
   + \sum_{t=1}^T  \underbrace{\vnorm{\mean{\v{Z}}(t) - \mean{\v{Z}}(t-1)}^2}_{:=\mathscr{F}(t)}. 
\end{aligned}
\label{eq:PEN}
}
\end{align}
The inequality follows from the fact that for any two scalars $a,b$, we have
\vspace{-.05in}
{
\begin{gather}
\begin{gathered}
    a^2 - b^2 = 2b(a-b) + (a-b)^2,\\
    \left(\proj_{\Rset_+}[a]\right)^2 - \left(\proj_{\Rset_+}[b]\right)^2
    \leq 2\proj_{\Rset_+}[b] (a-b) + (a-b)^2.
\end{gathered}
\end{gather}}

We separately bound $\mathscr{E}(t)$ and $\mathscr{F}(t)$.
For the former, we use the convexity of $\v{g}_j$ and the $\v{x}$-update to infer
\vspace{-.05in}
\begin{align}
{
\begin{aligned}
\mean{\v{Z}}(t) - \mean{\v{Z}}(t-1)
& = \frac{t}{N} \sum_{j=1}^N  \v{g}_j(\v{x}_j(t))
    - \frac{t-1}{N} \sum_{j=1}^N  \v{g}_j(\v{x}_j(t-1)) 
\\
& = \frac{t}{N} \sum_{j=1}^N  \v{g}_j\left( \frac{t-1}{t}\v{x}_j(t-1) + \frac{1}{t} \v{X}_j(t) \right)
\\
& \qquad - \frac{t-1}{N} \sum_{j=1}^N  \v{g}_j(\v{x}_j(t-1)) 
\\
& \leq 
\frac{1}{N} \sum_{j=1}^N  \v{g}_j\left(  \v{X}_j(t) \right).
\end{aligned}
\label{eq:Escr.bound}
}
\end{align}
Note that if an entry of $\v{g}$ encodes an equality constraint, the linearity of that constraint makes the above relation being met with an equality. Thus, we obtain
\vspace{-.05in}
\begin{align}
{
\begin{aligned}
    \mathscr{E}(t)
    &\leq
      \frac{1}{N}  \sum_{j=1}^N  \v{g}_j\left(  \v{X}_j(t) \right)^\T \proj_\Zset \left[ \mean{\v{Z}}(t-1)  \right].
\end{aligned}
}
\end{align}
To bound $\mathscr{F}_t$, we use the first line of \eqref{eq:Escr.bound} and the bounded/Lipschitz nature of $\v{g}_j$ on $\Xset_j$ to get
\vspace{-.05in}
\begin{align}
{
\begin{aligned}
    \vnorm{\mean{\v{Z}}(t) - \mean{\v{Z}}(t-1)}
    & \leq \frac{L_g}{N}\sum_{j=1}^N  (t-1) \vnorm{ \v{x}_j(t) - \v{x}_j(t-1)}
    + D_g
    \\
    & = \frac{L_g}{N}\sum_{j=1}^N   \vnorm{ \v{X}_j(t) - \v{x}_j(t)}
    +  D_g
    \\
    & \leq L_g  D_\Xset +  D_g
    \\
    & = D_Z.
\end{aligned}
\label{eq:Zbart.t1}
}
\end{align}
Replacing the bounds on $\mathscr{E}(t)$ and $\mathscr{F}(t)$ in \eqref{eq:PEN} gives the required bound on constraint violation in \eqref{eq:const.viol}, completing the proof of step (b).

\noindent $\bullet$ \emph{Step (c): Bounding the deviation of $\v{Z}_j$'s from its centroid:}
Consider $\v{\zeta} \in \Rset^{N \times M}$, given by
\vspace{-.05in}
\begin{align}
{
    \v{\zeta}(t)^\T = 
    \begin{pmatrix}
    \v{Z}_1(t) \mid \ldots \mid \v{Z}_N(t)
    \end{pmatrix}.
    \label{eq:zeta.def}
}
\end{align}
and define $\v{\Delta} := \v{I} - \frac{1}{N}\bone \bone^\T$, where $\bone \in \Rset^{N}$ is a vector of all ones and $\v{I}\in\Rset^{N \times N}$ is the identity matrix. Using this notation, we deduce
\vspace{-.05in}
\begin{align}
{
\sum_{j=1}^N \vnorm{ \v{Z}_j(t) -  \mean{\v{Z}}(t)}_2
 \leq \sqrt{N} \vnorm{\v{\Delta} \v{\zeta}(t)}_F  
\leq  N \vnorm{\v{\Delta} \v{\zeta}(t)}_2, 
\label{eq:norm2}
}
\end{align}
where, $\vnorm{\cdot}_F$ denotes the Frobenius norm of a matrix. Then, the $\v{Z}$-updates can be written as
\begin{align}
{
\v{\zeta}(t+1) = \v{W}\v{\zeta}(t) + \v{\varphi}(t), \quad \v{\zeta}(0) = 0
\label{eq:zeta.phi}
}
\end{align}
with $\v{\varphi}(t) \in \Rset^{N \times M}$; an analysis similar to \eqref{eq:Zbart.t1} gives that each row has a 2-norm bounded above by $D_Z$, implying
\begin{align}
{
    \begin{aligned}
    \vnorm{\v{\varphi}(t)}_2 \leq \sqrt{N} D_Z.
\end{aligned}
\label{eq:phi.bound}
}
\end{align}
Using \eqref{eq:zeta.phi}, we then obtain
\vspace{-.05in}
\begin{align}
{
\begin{aligned}
\vnorm{\v{\Delta} \v{\zeta}(t+1)}_2
= \vnorm{ \v{\Delta} \left( \v{W} \v{\zeta}(t)  + \v{\varphi}(t) \right)}_2 
\leq \vnorm{ \v{W} \v{\Delta} \v{\zeta}(t) }_2 + \vnorm{ \v{\Delta} \v{\varphi}(t)}_2
\end{aligned}
\label{eq:PYt1}
}
\end{align}
utilizing the fact that $\v{W}$ and $\v{\Delta}$  commute. To bound the first term in \eqref{eq:PYt1}, note that 
$\v{W}$ is doubly stochastic for which the Perron-Frobenius theorem \cite[Theorem 8.4.4]{HJbook} implies that its eigenvalue with largest absolute value is unity for which $\bone$ is the eigenvector. However, $\bone^\T \v{\Delta} = 0$, which in turn suggests $\v{\Delta} \v{\zeta}(t)$ is orthogonal to this eigenvector. Using the Courant-Fischer theorem \cite[Theorem 4.2.11]{HJbook}, we then obtain
\vspace{-.05in}
\begin{align}
{
    \vnorm{ \v{W} \v{\Delta}  \v{\zeta}(t) }_2
    \leq  \sigma_2(\v{W}) \vnorm{\v{\Delta}  \v{\zeta}(t) }_2  \label{eq:Xk1.2},
}
\end{align}
where $\sigma_2(\v{W})$ is the second largest singular value of $\v{W}$. Since $\v{W}$ is irreducible and aperiodic, $\sigma_2(\v{W}) \in (0,1)$.
We bound the second term in \eqref{eq:PYt1} as
\vspace{-.05in}
\begin{align}
{
\vnorm{ \v{\Delta} \v{\varphi}(t) }_2 
\leq \underbrace{\vnorm{\v{\Delta}}_2}_{= 1} \vnorm{\v{\varphi}(t)}_2
\leq \sqrt{N}D_Z,
\label{eq:VX}
}
\end{align}
because the 2-norm is sub-multiplicative.
Using the bounds in \eqref{eq:Xk1.2} and \eqref{eq:VX} in \eqref{eq:PYt1}, imply
\vspace{-.05in}
\begin{align}
{
\vnorm{\v{\Delta}\v{\zeta}(t+1)}_2
\leq \sigma_2(\v{W}) \vnorm{ \v{\Delta} \v{\zeta}(t)}_2
+ \sqrt{N} D_Z.
}
\end{align}
Iterating the above inequality gives
\vspace{-.05in}
\begin{align}
{
\begin{aligned}
\vnorm{\v{\Delta}\v{\zeta}(t)}_2
\leq 
 \sqrt{N} D_Z \sum_{\ell=0}^{t-1}  [\sigma_2(\v{W})]^{t-\ell-1}
\leq \sqrt{N} D_Z \left( 1-\sigma_2(\v{W}) \right)^{-1}.
\end{aligned}
\label{eq:PYt.2}
}
\end{align}
Then, \eqref{eq:norm2} and \eqref{eq:PYt.2} imply \eqref{eq:Z.centroid}, finishing step (c) of the proof.


$\bullet$ \emph{Step (d). Combining steps (a), (b), (c) to derive the result:}
Note that \eqref{eq:opt.gap} and \eqref{eq:const.viol} together with the definition of $\mathscr{V}_T$ give
\vspace{-.05in}
\begin{align}
{
\begin{aligned}
        \mathscr{V}_{T}(\v{x}(T), \v{z}(T))
        & 
        ={\sum_{j=1}^N \left[ f_j(\v{x}_j(T))  - \Dcal_j(\mean{\v{z}}(T)) \right]}
        + \frac{\eta T}{2N} \vnorm{\proj_\Zset \left[ \sum_{j=1}^N \v{g}_j(\v{x}_j(T)) \right]}^2
        \\
        &\leq \frac{2D_g}{T} \sum_{j=1}^N 
        \sum_{t=1}^T   \eta \vnorm{ \v{Z}_j(t-1) -  \mean{\v{Z}}(t-1)}
        + \frac{\eta}{2} N D_Z^2
        \\
        &\leq \eta \frac{ 2D_g
         N^{3/2} D_Z}{1-\sigma_2(\v{W})}
        + \frac{\eta}{2} N D_Z^2.
\end{aligned}
\label{eq:opt.gap.const.viol.Z.centroid}
}
\end{align}
where the second inequality follows from using \eqref{eq:Z.centroid}. Using $\eta = \eta_0/\sqrt{T}$, we then obtain the upper bound in \eqref{eq:result}.

\subsection{Lower bounding $\mathscr{V}_T$}
By the saddle-point property of a primal dual optimizer $(\v{x}^\star, \v{z}^\star)$ of $\Pcal$, we get
\vspace{-.05in}
\begin{align}
    {
    \begin{aligned}
    \Pcal^\star 
    &= \Lcal(\v{x}^\star, \v{z}^\star)
    \\
    &\leq \Lcal(\v{x}(T),\v{z}^\star)
    \\
    &= \sum_{j=1}^N f_j(\v{x}_j(T)) + \v{z}^{\star, \T} \sum_{j=1}^N \v{g}_j(\v{x}_j(T))
    \\
    &\leq \sum_{j=1}^N f_j(\v{x}_j(T)) + \v{z}^{\star, \T} \proj_{\Zset}\left[ \sum_{j=1}^N \v{g}_j(\v{x}_j(T))\right].
    \end{aligned}
    }
\end{align}
Applying Young's inequality to the last summand in the right hand side of the above relation, we further get
\vspace{-.05in}
\begin{align}
{
    \begin{aligned}
    \Pcal^\star 
    \leq \sum_{j=1}^N f_j(\v{x}_j(T)) + \frac{N}{2 \eta T}\vnorm{\v{z}^\star}^2 
    + \frac{\eta T}{2N} \vnorm{\proj_{\Zset}\left[ \sum_{j=1}^N \v{g}_j(\v{x}_j(T))\right]}^2.
    \end{aligned}
}
\end{align}
Subtracting $\sum_{j=1}^N \Dcal_j(\mean{\v{z}}(T))$ on both sides and using $\eta = \eta_0/\sqrt{T}$ yields the desired lower bound on $\mathscr{V}_T$ in \eqref{eq:result}.

To summarize, the lower and upper bound of $\mathscr{V}_{T}(\v{x}(T), \mean{\v{z}}(T))$ are
\begin{align}
    \Pcal^\star - \sum_{j=1}^N \Dcal_j(\mean{\v{z}}(T))  - \frac{N}{2 \eta T}\vnorm{\v{z}^\star}^2 
    \leq \mathscr{V}_{T}
    \leq \eta \frac{ 2D_g N^{3/2} D_Z}{1-\sigma_2(\v{W})} + \frac{\eta}{2} N D_Z^2.
\end{align}
with $\eta = \eta_0/\sqrt{T}$. Note $D_g$, $D_Z$ and $N$ reflect the required iterations that are affected by the problem structure.
One expects $D_g$, and as a result, $D_Z$ to grow with the number of constraints \eqref{eq:P.E}, \eqref{eq:P.I}, and $N$ to grow with the number of agents.

%% file: multi_area.tex

\section{Grid Optimization Example 1: Tie-Line Scheduling}
\label{sec:multi_area}
In this section and the next two, we present three different examples of grid optimization problems that can be cast as examples of $\Pcal$. { See \cite{extended} for a literature survey on each problem. Here, we focus instead on reformulation of each as an instance of $\Pcal$ and report numerical results from running Algorithm \ref{alg:ddsa} to reveal interesting properties of the algorithm.}

We first present our results on $\Pcal_1$--the tie-line scheduling problem.
Tie-lines are transmission lines that interconnect the footprints of different system operators (henceforth, called areas). 
Ideally, one would solve a joint OPF problem over assets within all areas to optimize tie-line schedules. However, technical and legal challenges impede aggregation of all problem data at a central location, requiring a distributed algorithm design.
{ Since the seminal work in \cite{KimB1997}, a substantial literature has developed on multi-area OPF problems for tie-line scheduling; see \cite{extended} for a survey of methods. Here, we formulate $\Pcal_1$ and apply Algorithm \ref{alg:ddsa}.
}
{Denote by 
$\v{\theta}_j \in \Rset^{n_j}$ and $\v{\theta}_{\overline{j}} \in \Rset^{n_{\overline{j}}}$, the voltage phase angles at the internal and the boundary buses in each area $j$, respectively. We adopt a linear power flow model in which the vector of power injections within an area (generation less demand $\v{p}^G_j - \v{p}^D_j$ at internal buses and zero injections at boundary buses) become linear in voltage phase angles through suitably defined matrices $\v{B}$. Individual line flows within and across areas also become linear in these angles, defined through matrices $\v{H}$. Angles are constrained within $[0, 2\pi]$, represented as $\v{\theta} \in \Theta$.
Utilizing $\v{L}$'s to encode line capacities, $c$'s to denote power procurement costs (typically deduced from supply offers in electricity markets), and $\Gfk(N, \Eset)$ to represent the interconnection graph among the areas, the multi-area OPF problem becomes}
\begin{subequations}
\begin{alignat}{2}
\Pcal_1 : \ & \text{minimize} \ \ &&
	  \sum_{j=1}^N c_j\left(\v{p}_j^G \right), \notag\\
 & \text{subject to} && \underline{\v{p}}_j^G \leq \v{p}_j^G \leq \overline{\v{p}}_j^G, \v{\theta} \in \Theta,
 \label{eq:MA.1}
 \\
&&& \v{B}_{j,j} \v{\theta}_j + \v{B}_{j,\overline{j}} \v{\theta}_{\overline{j}} = \v{p}_j^G - \v{p}_j^D, 
 \label{eq:MA.2}
 \\
&&& \v{B}_{\overline{j},j} \v{\theta}_j +  \v{B}_{\overline{j},\overline{j}} \v{\theta}_{\overline{j}} +  \sum_{k \sim j}\v{B}_{\overline{j},\overline{k}} \v{\theta}_{\overline{k}} = 0,
 \label{eq:MA.3}
 \\
&&& \v{H}_{j} \v{\theta}_j + \v{H}_{ \overline{j}} \v{\theta}_{\overline{j}}  \leq \v{L}_j, 
\label{eq:MA.4}
 \\
&&& \v{H}_{j,k} \v{\theta}_{\overline{j}} + \v{H}_{k,j} \v{\theta}_{\overline{k}}  \leq \v{L}_{jk}, 
\label{eq:MA.5}
\\
&&&  j=1,\ldots,N, \ k \sim j \text{ in } \Gfk. \notag
\end{alignat}
\label{eq:MAOPF}
\end{subequations}
Here, \eqref{eq:MA.1}--\eqref{eq:MA.4} encode the generation capacity and angle constraints, power balance and transmission line constraints within each area, while \eqref{eq:MA.5} enforces limits on  tie-line flows.
To cast \eqref{eq:MAOPF} as $\Pcal$, define
\begin{gather*}
\v{x}_{j} 
= \left( \v{\theta}_{j}^\T, \v{\theta}_{\overline{j}}^\T, [\v{p}_j^G]^\T \right)^\T, 
\
\Xset_j 
= \left\{ \v{x}_j \ | \ \eqref{eq:MA.1}, \eqref{eq:MA.2}, \eqref{eq:MA.4} \right\}, 
\ 
f_j(\v{x}_j) 
= c_j\left(\v{p}_j^G\right),
\end{gather*}
and write \eqref{eq:MA.3} and \eqref{eq:MA.5} as \eqref{eq:P.E} and \eqref{eq:P.I}, respectively.

\begin{figure}[hbtp]
    \centering
        \includegraphics[width=0.6\textwidth]{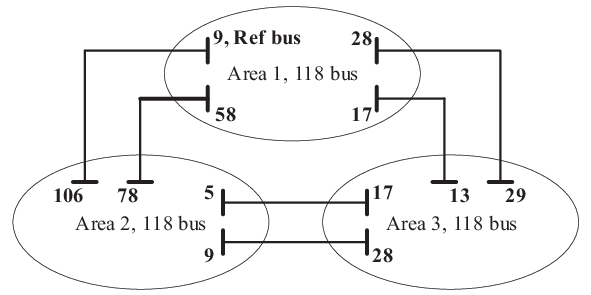}
        \vspace{-.05in}
        \caption{The three-area network for multi-area optimal power flow simulations, obtained by joining three IEEE 118-bus systems.}
		\label{fig:ma.1}    
    \vspace{.1in}
	\includegraphics[width=0.8\textwidth]{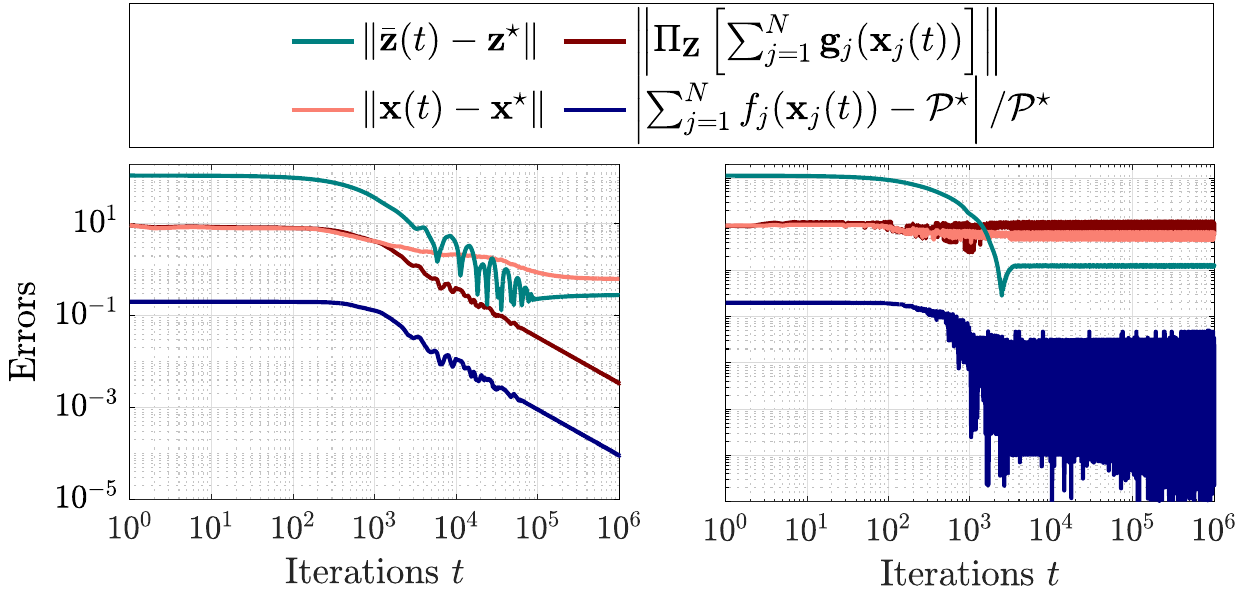}
	\vspace{-.05in}
	\caption{
	Performance of Algorithm \ref{alg:ddsa} (left) and Algorithm \ref{alg:dualsubg} (right) on $\Pcal_1$ for the network in Figure \ref{fig:ma.1}.}	
	\label{fig:ma.2}
	    \vspace{0.1in}
	    \includegraphics[width=0.8\textwidth]{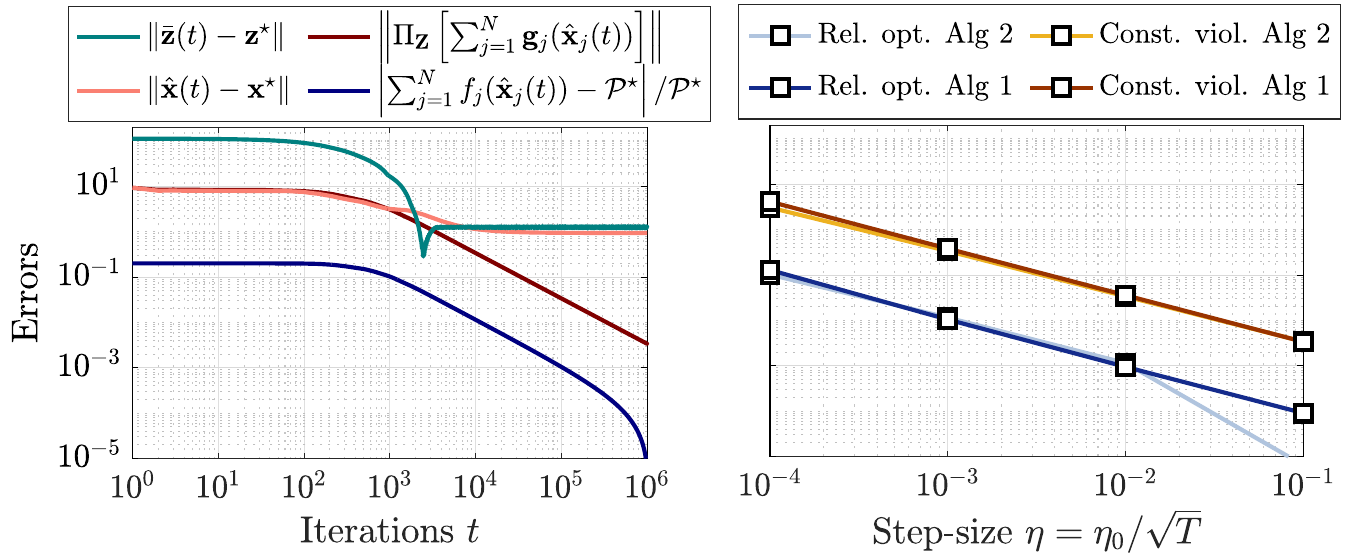}
	\caption{Performance of Algorithm \ref{alg:dualsubg} with primal averaging (left) and the impact of step-size on  Algorithms \ref{alg:ddsa} and \ref{alg:dualsubg} with primal averaging (right).}
	\label{fig:ma.3}
	\vspace{-0.1in}
\end{figure}

Consider the three-area power system shown in Figure \ref{fig:ma.1} that comprises three IEEE 118 systems stitched together with 6 tie-lines as shown. The three systems were modified as delineated in \cite[Appendix B1]{extended}. We applied Algorithm \ref{alg:ddsa} on a reformulation of $\Pcal_1$ as an instance of $\Pcal$ with a flat start ($\v{z}_j(1) = 0$, $\v{Z}_j(0) = 0$, $j = 1, \ldots, N$) and step size $\eta = \eta_0/\sqrt{T}$, where $\eta_0 = 10^{2}$ and $T=10^{6}$. The results are portrayed on the left of Figure \ref{fig:ma.2}. We chose 
$\v{W}$ based on the transition probabilities of a Markov chain in the Metropolis-Hastings algorithm (see \cite[Sec. 2.5]{notarstefanoDistributedOptimizationSmart2020}). Here, $\Pcal^\star$ was computed by solving $\Pcal_1$ as a linear program. Our simulations were performed in MATLAB 2018b.
All sub-problems were solved via MOSEK version 9.2.17.

\begin{algorithm}
	\caption{Distributed dual subgradient to solve $\Pcal$.}
	\label{alg:dualsubg}

    Choose $\v{z}_j(1) = 0$ and $\eta = \eta_0/\sqrt{T}$.

	\For{$t=1, \ldots, T$}{
	
	$\v{x}_j(t)
    \gets \argmin_{\v{x}_j \in \Xset_j} \ \Lcal_j(\v{x}_j, \v{z}_j(t))$. \label{alg2:x.2}
    
    $\v{z}_j({t+1})
    \gets \sum_{k=1}^N W_{jk} \proj_\Zset\left[ \v{z}_k(t) + \eta \v{g}_k(\v{x}_k(t))\right]$. \label{alg2:z.2}
    }

\end{algorithm}

We compared Algorithm \ref{alg:ddsa} with the classical dual subgradient method in Algorithm \ref{alg:dualsubg} (the projection and the consensus operations in step \ref{alg2:z.2} are sometimes reversed, e.g., in \cite{simonetto2016primal}). The progress of Algorithm \ref{alg:dualsubg} with the same step-size used for Algorithm \ref{alg:ddsa} are shown in the right of Figure \ref{fig:ma.2}. Note that Algorithm \ref{alg:ddsa} leads to much smoother progress of $\sum_{j=1}^N f_j(\v{x}_j(t))$ compared to that with Algorithm \ref{alg:dualsubg}. Classical dual subgradient with primal averaging via
$\hat{\v{x}}_j(t) := \frac{1}{t}\sum_{r=1}^t \v{x}_j(r)$ for each $j = 1, \ldots, N$ can prevent this ``flutter'' (see \cite[Section 4]{NedicO2009_PD}), as the left plot in Figure \ref{fig:ma.3} reveals. While step \ref{alg:x} of Algorithm \ref{alg:ddsa} executes a similar averaging operation, this averaging step cannot be viewed as an \emph{output} of the iteration dynamics as is the case for Algorithm \ref{alg:dualsubg} with averaging. As a result, the last iterate of Algorithm \ref{alg:ddsa} moves smoothly as opposed to Algorithm \ref{alg:dualsubg}. Such an update is useful in applications that require iterates to be directly implemented as control actions and the dual subgradient is only available at the current iterate (see \cite{nesterov2018dual} for a detailed discussion).  


In the right of Figure \ref{fig:ma.3}, we compared the impact of step-size on the performance of Algorithms \ref{alg:ddsa} and  \ref{alg:dualsubg} with primal averaging after $T=10^{6}$ iterations. 
Here, relative optimality measures $\left | \sum_{j=1}^N f_j(\v{x}_j(t)) - \Pcal^\star \right | / \Pcal^\star $ and constraint violation measures $\vnorm{\proj_\Zset[\sum_{i=1}^N\v{g}_j(\v{x}_j(T))]}$ for Algorithm \ref{alg:ddsa} and the same at $\hat{\v{x}}_j(t)$'s for Algorithm \ref{alg:dualsubg}. Empirically, constraint violation for Algorithm \ref{alg:ddsa} appears similar to that for Algorithm \ref{alg:dualsubg} with primal averaging that is known to have $\Ocal(T^{-1/2})$ decay rate, much better than that suggested by \eqref{eq:result.const.viol} for Algorithm \ref{alg:ddsa}.

%% file: DER_coordination.tex
\section{Grid Optimization Example 2: DER Coordination}
\label{sec:DER_coord}
Our next application problem is the coordination of DERs such as thermostatically controlled loads, electric vehicles, distributed rooftop solar, etc. that are increasingly getting adopted in distribution grids. There is a long literature on DER coordination to fulfill a variety of objectives that range from tracking a regulation signal at the T\&D interface, to volt/VAR control within the distribution grid, etc. (e.g., see \cite{wuHierarchicalControlFramework2017,DallAnese2018})
{ See \cite{extended} for a longer list.}

We formulate the DER coordination problem $\Pcal_2$ over a balanced three-phase radial distribution network on $N$ buses, described by graph $\Gfk(N, \Eset)$. Let the first bus be the T\&D interface. Associate directions to edges in $\Eset$ arbitrarily to obtain a directed graph $\vec{\Gfk}(N, \vec{\Eset})$, where $j \to k \in \vec{\Eset}$ denotes a directed edge from bus $j$ to bus $k$ in $\vec{\Gfk}$. At each bus $j$, consider a dispatchable asset capable of injecting real and reactive powers $p^G_j, q^G_j$, respectively. Let $c_j(p^G_j, q^G_j)$ denote the cost of power procurement from that dispatchable asset. At $j=1$, this cost might reflect the cost of procuring power from the transmission grid. It can also encode deviation of the power injection from a set point defined by a frequency regulation signal. At buses $j=2,\ldots,N$, the cost can encode the disutility of deferred demand or cost of power production from dispatchable generation. The power injection capabilities of this asset at bus $j$ are limited as $\underline{p}_j^G \leq p_j^G \leq \overline{p}_j^G$ along with 
\begin{gather*}
\underline{q}_j^G \leq q_j^G \leq \overline{q}_j^G \text{ or } 
 \left[ p_j^G \right]^2 + \left[ p_j^G \right]^2  \leq \left[\overline{s}_j^G\right]^2,
 \end{gather*}
henceforth denoted as $(p_j^G, q_j^G) \in \Sset_j$.
Such models encompass photovoltaic and energy storage systems, water pumps, commercial HVAC systems, etc. At each bus $j$, also assume nominal real and reactive power demands $p^D_j$ and $ q^D_j$.

We need additional notation to describe the DER coordination problem. Associate with bus $j$ the squared voltage magnitude $w_j$, { deemed to lie in $[\underline{w}_j, \overline{w}_j]$.} Let $P_{j,k}$, $Q_{j,k}$ denote the real and reactive power flows from bus $j$ to bus $k$ for $j\to k$ in $\vec{\Gfk}$. Denote by $\ell_{j,k}$, the squared current magnitude flowing from bus $j$ to bus $k$, { upper bounded by $L_{j, k}$}.  
Let $\mathpzc{r}_{j,k}$ and $\mathpzc{x}_{j,k}$ denote the resistance and reactance of the line $j \to k$. 
The DER coordination problem with a second-order conic  relaxation of power flow equations in the radial distribution network can be formulated as
\vspace{-.05in}
\begin{subequations}
\begin{alignat}{2}
\Pcal_2: \ & \text{minimize} \quad \sum_{j=1}^N  c_j(p^G_j, q^G_j), 
\notag
\\
& \text{subject to} \notag \\
& \qquad \left( p_{j}^G, q_{j}^G \right) \in \Sset_j, 
\label{eq:DER.1}
\\
& \qquad p_j^G - p_j^D = \sum_{k:j \to k}{P_{j,k}} - \sum_{k:k \to j}{(P_{k,j} - \mathpzc{r}_{k,j}\ell_{k,j})}, 
\label{eq:DER.2}
 \\
& \qquad q_j^G - q_j^D = \sum_{k:j \to k}{Q_{j,k}} - \sum_{k:k \to j}{(Q_{j,k} - \mathpzc{x}_{k,j} \ell_{k,j})}, 
\label{eq:DER.3}
\\
& \qquad w_k = w_j  - 2(\mathpzc{r}_{j,k} P_{j,k} + \mathpzc{x}_{j,k} Q_{j,k})
+ (\mathpzc{r}_{j,k}^2 + \mathpzc{x}_{j,k}^2) \ell_{j,k}, 
\label{eq:DER.4}
\\
& \qquad \ell_{j,k} \leq L_{j,k},\ \underline{w}_j \leq w_j \leq \overline{w}_j, 
\label{eq:DER.5}
\\
& \qquad \ell_{j,k} w_j  \geq P_{j,k}^2 + Q_{j,k}^2, 
\label{eq:DER.6}
\\
& \qquad j = 1,\ldots, N, j\to k \in \ \vec{\Gfk}. \notag
\end{alignat}
\label{eq:DER_coord}
\end{subequations}
The last inequality is a second-order cone constraint, making \eqref{eq:DER_coord} a second-order cone program (SOCP). 
Constraints in \eqref{eq:DER.2}, \eqref{eq:DER.3} describe real and reactive power balance at each bus. Relations \eqref{eq:DER.2}, \eqref{eq:DER.3}, \eqref{eq:DER.3}, and the inequality in \eqref{eq:DER.6} replaced with an equality define the feasible set described by AC power flow equations (see \cite{Farivar2013,low2014convex} for details). 

To cast $\Pcal_3$ as an instance of $\Pcal$, we first write the out-neighbors of $j$ in $\vec{\Gfk}$ as $k_1, \ldots, k_J$ and identify 
\begin{gather*}
\v{x}_j := \left( p_j^G, q_j^G, w_j, P_{j, k_1}, \ldots, P_{j, k_{J}}, 
 Q_{j, k_1}, \ldots, Q_{j, k_J}, \ell_{j, k_1}, \ldots, \ell_{j, k_{J}} \right)^\T, \\
\Xset_j := \{ \v{x}_j \ | \ \eqref{eq:DER.1}, \eqref{eq:DER.5}, \eqref{eq:DER.6} \}, \
f_j(\v{x}_j) = c_j(p_j^G, q_j^G).
\end{gather*}
Then, it is straightforward to write \eqref{eq:DER.2}, \eqref{eq:DER.3} and \eqref{eq:DER.4} as examples of \eqref{eq:P.E}. This formulation does not require inequality constraints of the form \eqref{eq:P.I}.
Note T\&D interface's energy balancing and voltage constraints \eqref{eq:DER.2}-\eqref{eq:DER.4} can be treated as \eqref{eq:P.X} since there is no coupling variable from adjacent agents.


We ran Algorithm \ref{alg:ddsa} on $\Pcal_2$ over a modified IEEE 4-bus radial distribution network (see \cite[Appendix B2]{extended} for details). 
\begin{figure}[htbp]
    \centering 
    \vspace{-0.15in}
    \subfloat[]{
        \includegraphics[width=0.2\textwidth]{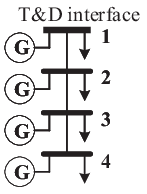}  \label{fig:der.1}
    }  
    \subfloat[]{
        \includegraphics[width=0.4\textwidth]{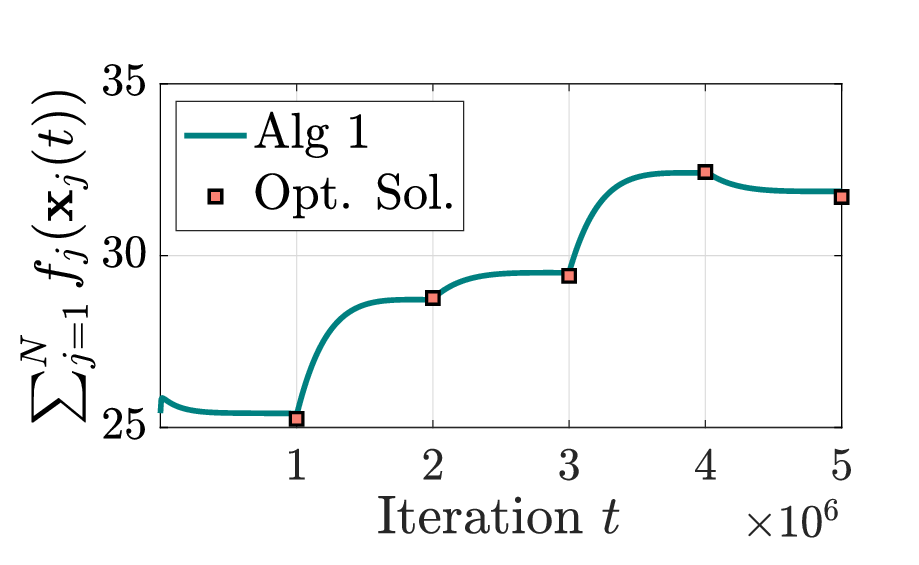}     
        \label{fig:der.2}
    }        
    \caption{
    (a) A 4-bus radial network.
    (b) Progress of the objective function at the last iterate of Algorithm \ref{alg:ddsa} on $\Pcal_2$ for the network in Figure \ref{fig:der.1}.
    }
    \vspace{-0.1in}
\end{figure}


System conditions in the distribution grid can change quite fast. One line of work on DER coordination solves optimization problems in quick successions to deal with such changes, e.g., in \cite{Zhou2017}. 
To illustrate the use of DER coordination with time-varying distribution grid conditions, we simulated a case where real and reactive power demands were changed every $10^{6}$ iterations as prescribed in \cite[Appendix B2]{extended} with step-size $\eta = 0.1$.
Algorithm \ref{alg:ddsa} is restarted after every change. 
Here, we use the last primal-dual iterate at the point of change to restart Algorithm \ref{alg:ddsa}. As Figure \ref{fig:der.2} illustrates, Algorithm \ref{alg:ddsa} can track the optimal cost in the changing problem environment.

%



\begin{figure}[htbp]
    \centering   
        \vspace{-.05in}
        \includegraphics[width=0.8\textwidth]{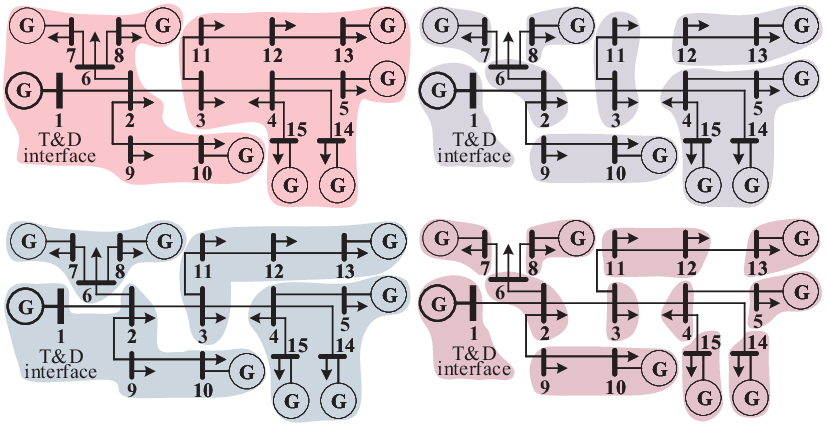}
        \caption{The IEEE 15-bus test feeder subdivided into 2, 4, 8 and 12 groups.}
    \label{fig:der.networks}
        \includegraphics[width=0.7\textwidth]{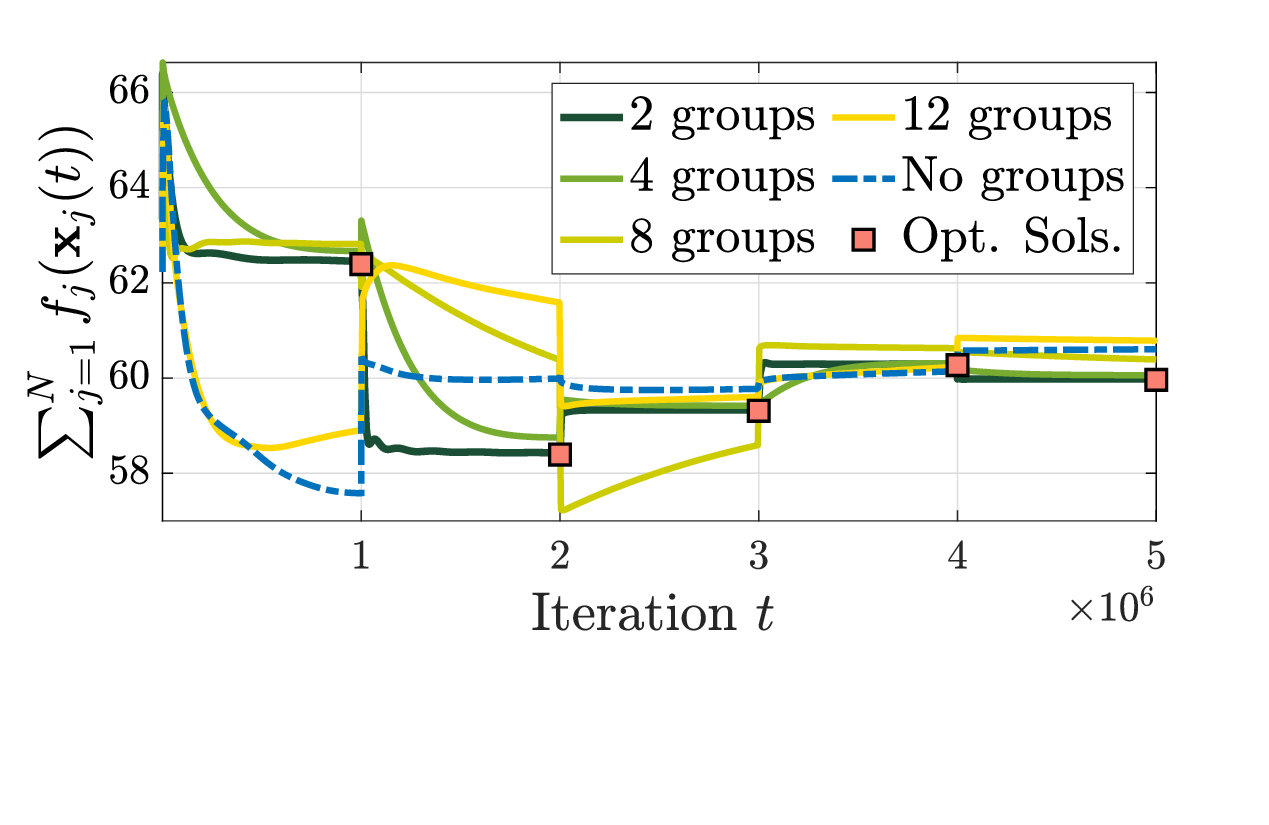} 
        \vspace{-30pt}
    \caption{Evolution of the objective function value of Algorithm \ref{alg:ddsa} on the IEEE 15-bus test system with varying degrees of decentralization (based on groupings of buses per Figure \ref{fig:der.networks}).} 
    \label{fig:der.5}
\end{figure}

{
Convergence slows down over larger networks. Consider a modified IEEE 15-bus radial distribution network; see \cite[Appendix B2]{extended} for the modifications. Figure \ref{fig:der.5} illustrates that Algorithm \ref{alg:ddsa} with restarts is able to track optimal costs, but only when nodes are grouped together into fewer agents in  $\Pcal$, as shown in Figure \ref{fig:der.networks}. Such a slowdown is expected, given that the convergence guarantees depend on $N$, albeit polynomially. 
We remark that \emph{all} first-order algorithms suffer from this issue. In practical implementation, one must carefully explore the trade-off between the degree of decentralization and the accuracy of tracking within a fixed number of iterations.
}

%% file: TD_coordination.tex
\section{Grid Optimization Example 3: T\&D Coordination}
\label{sec:TD_coordination}

The lack of visibility of transmission SOs into distribution grids and bottlenecks in wholesale market clearing software make it impossible for such SOs to directly harness the flexibility offered by DERs in the distribution networks.
Naturally, distributed algorithms are suited for T\&D coordination; see \cite{extended} for prior art.
Assume that for each distribution network, an aggregator $\Acal$ directly controls the dispatchable DERs and knows the network parameters. In what follows, we present the T\&D coordination problem $\Pcal_3$ that a transmission SO and a collection of distribution aggregators solve in a distributed fashion. We utilize a semidefinite relaxation of power flow equations for the transmission network and a linear distribution flow model for the distribution grids. 


\newcommand{\tran}{{\textrm{tran}}}
\renewcommand{\dist}{{\textrm{dist}}}

To formulate the joint dispatch problem of all T\&D assets, we require three different graphs. The first among these is the transmission network, modeled as an undirected graph $\Gfk^\tran$ on $n^\tran$ transmission buses. The second set of graphs are the distribution grids that connect to the transmission network at their points of common coupling--the $n^\tran$ transmission buses. 
We model the distribution grid connected to transmission bus $\ell$ as an undirected graph $\Gfk^\dist_\ell$ on $n^\dist_\ell + 1$ distribution buses, where the first bus of $\Gfk^\dist_\ell$ coincides with bus $\ell$ in $\Gfk^\tran$. Finally, we consider an undirected star graph $\Gfk$ on $N = n^\tran+1$ nodes with
the aggregators $\Acal_1, \ldots, \Acal_{n^\tran}$ as the satellite nodes and the SO (the $N$-th node) at the center.

Let $\v{V} \in \Cset^n$ denote the vector of nodal voltage phasors, where $\Cset$ is the set of complex numbers. We formulate the engineering constraints of the grid using the positive semidefinite matrix $\v{\Lambda} := \v{V} \v{V}^{\Hsf} \in \Cset^{n^\tran \times n^\tran}$. To describe these constraints, let $y_{\ell,k} = y_{k,\ell}$ denote the admittance of the transmission line joining buses $\ell,k$ in $\Gfk^\tran$ and $y_{\ell, \ell}$ denote the shunt admittance at bus $\ell$. Then, define $\v{\Phi}_{\ell,k}$, $\v{\Psi}_{\ell,k}$ as the $n^\tran \times n^\tran$ Hermitian matrices whose only nonzero entries are 
\begin{gather*}
[\v{\Phi}_{\ell, k}]_{\ell, \ell} := \frac{1}{2}(y_{\ell,k} + y_{\ell,k}^{\Hsf}),  
[\v{\Phi}_{\ell, k}]_{\ell, k} =  [\v{\Phi}_{\ell,k}]_{k,\ell}^{\Hsf} := -\frac{1}{2}y_{\ell,k}, \\
[\v{\Psi}_{\ell, k}]_{\ell, \ell} := \frac{1}{2 \ii}(y_{\ell,k}^{\Hsf} - y_{\ell,k}),  
[\v{\Psi}_{\ell, k}]_{\ell, k} =  [\v{\Psi}_{\ell,k}]_{k,\ell}^{\Hsf} := \frac{1}{2\ii}y_{\ell,k}.
\end{gather*}
In addition, we define the ${n^\tran \times n^\tran}$ Hermitian matrices 
\begin{align*}
\v{\Phi}_\ell
:=
\frac{1}{2}\left( {y_{\ell,\ell} + y_{\ell,\ell}^{\Hsf}}\right) \bone_\ell \bone_\ell^{\Hsf} 
+ \sum_{k \sim \ell} \v{\Phi}_{\ell,k}, 
\quad
\v{\Psi}_\ell 
:= 
\frac{1}{2\ii} \left( {y_{\ell,\ell}^{\Hsf} - y_{\ell,\ell}}\right) \bone_\ell \bone_\ell^{\Hsf}
+ \sum_{k \sim \ell} \v{\Psi}_{\ell,k},
\end{align*}
where $\bone$ is a vector of all ones of appropriate size and $\bone_\ell$ is a vector of all zeros except at the $\ell$-th position that is unity.
This notation allows us to describe the apparent power flow from bus $\ell$ to bus $k$ as $ \trace(\v{\Phi}_{\ell, k}  \v{\Lambda}) + \ii \trace( \v{\Psi}_{\ell, k}  \v{\Lambda} )$, the apparent power injection at bus $\ell$ as $\trace( \v{\Phi}_\ell \v{\Lambda}) + \ii \trace(\v{\Psi}_\ell \v{\Lambda})$,
and the squared voltage magnitude at bus $\ell$ as $\trace( \bone_\ell \bone_\ell^{\Hsf} \v{\Lambda})$, { where $\ii:=\sqrt{-1}$}. At each transmission bus $\ell$ with load $P_\ell^D + \ii Q_\ell^D$, let a generator supply apparent power $P_\ell^G + \ii Q_\ell^G$ with procurement cost described by $C_\ell$. 

Let each transmission bus $\ell$ be the first bus of an $n^\dist+1$-bus distribution network $\Gfk^\dist_\ell$. Let $\v{p}_\ell + \ii \v{q}_\ell \in \Cset^{n^\dist}$ denote the vector of net power injections across the distribution network, save the first bus. Further, let the power procurement cost be given by $\v{c}_\ell$ to inject $\v{p}_\ell + \ii \v{q}_\ell \in \Cset^{n^\dist}$.
Also, let $\v{\Lambda}_\ell \in \Cset^{n^\dist}$ denote the vector of squared voltage magnitudes across the same set of buses. We adopt the popular LinDistFlow model to tackle the nonconvex nature of the power flow equations in the distribution grid. Let $\widetilde{\v{M}} \in \Rset^{n^\dist \times n^\dist}$ be the node-to-edge incidence matrix of $\Gfk^\dist_\ell$. Further, remove the first row of $\widetilde{\v{M}}$ to obtain the reduced incidence matrix $\v{M}$. Then, the voltage magnitudes are related to power injections under the LinDistFlow model as
    $\v{\Lambda}_\ell = \v{\rho}_\ell \v{p}_\ell + \v{\chi}_\ell \v{q}_\ell + \Lambda_{\ell, \ell} \bone$,
where $\v{\rho}_\ell$ and $\v{\chi}_\ell$ are $n^\dist \times n^\dist$ matrices defined as 
$\v{\rho}_\ell := 2 \v{M}^{-\T} \diag({\v{\mathpzc{r}}_\ell}) \v{M}^{-1}$, $\v{\chi}_\ell := 2 \v{M}^{-\T} \diag({\v{\mathpzc{x}}_\ell}) \v{M}^{-1}$,
$\v{\mathpzc{r}}_\ell$/ $\v{\mathpzc{x}}_\ell$ collect the resistances/reactances of the $n^\dist$ distribution lines. 

The optimal joint dispatch over all T\&D assets is given by
\begin{subequations}
\begin{alignat}{2}
	\Pcal_3: \ & {\text{minimize}} && \quad  \sum_{\ell=1}^{n^\tran} C_{\ell} (P^G_{\ell}, Q^G_{\ell}) +   \sum_{\ell=1}^{n^\tran} \v{c}_{\ell}(\v{p}^G_{\ell}, \v{q}^G_{\ell}), 
	\notag
	\\
	& \text{subject to} 
	&& \quad ({P}_\ell^G, Q_\ell^G) \in \Sset^\tran_{\ell}, 
		\label{eq:TD.1}
		\\
	&&& \quad \left( \v{p}_{\ell}^G, \v{q}_{\ell}^G \right) \in \Sset^\dist_{\ell},
		\label{eq:TD.2}
		\\
	&&& \quad P^G_\ell - P_{\ell}^D + \bone^\T \left( \v{p}^G_{\ell} - \v{p}^D_{\ell} \right) = \trace(\v{\Phi}_\ell \v{\Lambda}),  
		\label{eq:TD.3}
		\\
	&&& \quad Q^G_\ell - Q_{\ell}^D + \bone^\T \left( \v{q}^G_{\ell} - \v{q}^D_{\ell} \right) = \trace( \v{\Psi}_\ell \v{\Lambda}),
		\label{eq:TD.4}
		\\
	&&& \quad \trace(\v{\Phi}_{\ell, \ell'} \v{\Lambda} ) \leq L_{\ell, \ell' }, 			
		\label{eq:TD.5}
		\\
	&&& \quad \ul{w}_\ell \leq  \Lambda_{\ell, \ell} \leq \ol{w}_\ell,  
		\label{eq:TD.6}
		\\			
	&&& \quad \v{\Lambda} \succeq 0,
		\label{eq:TD.7}
		\\
    &&& \quad  \ul{\v{w}}_\ell \leq          \v{\rho}_\ell \v{p}_\ell + \v{\chi}_\ell \v{q}_\ell + \Lambda_{\ell, \ell} \bone \leq \ol{\v{w}}_\ell,
        \label{eq:TD.8}     
        \\  
	&&& \quad \text{for } \ell = 1,\ldots,n^\tran, \ \ell' \sim \ell. \notag
\end{alignat}
\label{eq:TD}
\end{subequations}
Here, \eqref{eq:TD.1}, \eqref{eq:TD.2} encode the capabilities of transmission and distribution assets, respectively, where $\Sset^\tran$ and $\Sset^\dist$ are assumed convex. Constraints in  \eqref{eq:TD.3} and \eqref{eq:TD.4} enforce nodal power balance at transmission nodes. Transmission line (real power) flows are constrained in \eqref{eq:TD.5}, with $L$'s describing the line capacities. One can alternately constrain apparent power flows. Transmission voltage limits are enforced via \eqref{eq:TD.6} within $[\underline{w}_\ell, \overline{w}_\ell]$. The relation $\v{\Lambda}=\v{V} \v{V}^{\Hsf}$ requires $\v{\Lambda}$ to be positive semidefinite (enforced in \eqref{eq:TD.7}) and rank-1. We consider the rank-relaxed power flow equations by dropping the rank-1 requirement. For the distribution grid, voltage limits are enforced in \eqref{eq:TD.8} and power flow limits  are  ignored. 

Recall that $\Gfk$ for T\&D coordination problem is a graph on $N = n^\tran + 1$ nodes, where the first $n^\tran$ nodes are transmission buses and the last node represents the SO. Define
\begin{align*}
 \v{x}_\ell := \left( \v{p}_\ell^\T, \v{q}_\ell^\T \right)^\T, \
\Xset_\ell := \{ \v{x}_\ell \ | \  \eqref{eq:TD.2} \}, \
f_\ell = c_\ell(\v{p}^G_{\ell}, \v{q}^G_{\ell})
\end{align*} 
for $\ell = 1,\ldots,n^\tran$. Collect the real and reactive power generations across the transmission grid in the vectors $\v{P}^G, \v{Q}^G$, respectively. Then, define 
\begin{align*}
\v{x}_{N} &:= \left( [\v{P}^G]^\T, [\v{Q}^G]^\T, \vvec{\Re\{\v{\Lambda}\}}^\T, \vvec{\Im\{\v{\Lambda}\}}^\T \right)^\T, 
\\
\Xset_{N} &:= \{ \v{x}_{n^\tran + 1} \ | \ \eqref{eq:TD.1}, \eqref{eq:TD.5}, \eqref{eq:TD.6}, \eqref{eq:TD.7}\}, 
\\
f_{N}(\v{x}_{N}) &= \sum_{\ell=1}^{n^\tran} C_\ell (P_\ell^G, Q_\ell^G).
\end{align*}
The constraint \eqref{eq:TD.7} can be written in terms of $\v{x}_N$ as
$$ \begin{pmatrix} \Re\{ \v{\Lambda}\} & \Im\{ \v{\Lambda}\} \\ - \Im\{ \v{\Lambda}\} & \Re\{ \v{\Lambda}\}\end{pmatrix} \succeq 0$$
and \eqref{eq:TD.3} -- \eqref{eq:TD.4} as examples of \eqref{eq:P.E} using 
$\trace (\v{\varphi} \v{\Lambda} ) 
= \vvec{\Re\{\v{\varphi}\}}^\T \vvec{\Re\{\v{\Lambda}\}} 
+ 
\vvec{\Im\{\v{\varphi}\}}^\T \vvec{\Im\{\v{\Lambda}\}}$
for a Hermitian matrix $\v{\varphi}$. Constraints \eqref{eq:TD.8} are examples of inequality constraints in \eqref{eq:P.I}.

\begin{figure}[htbp]
    \centering 
    \vspace{-0.05in}
    \includegraphics[width=0.8\textwidth]{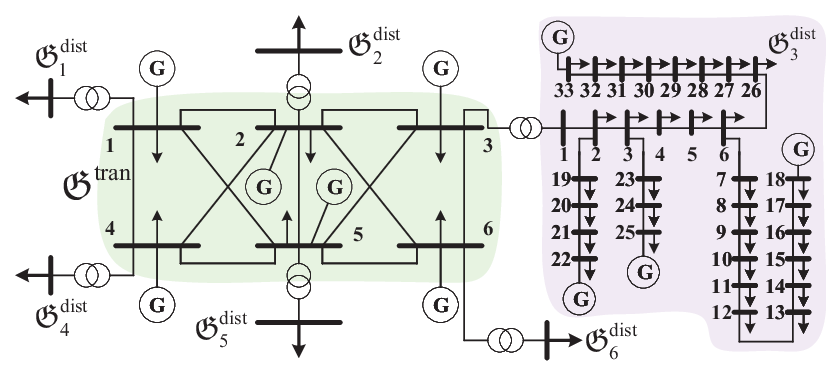}     
    \caption{The 204-bus network for T\&D simulations, obtained by joining the IEEE 6-bus transmission network with six IEEE 33-bus distribution networks.}
    \label{fig:td.1}
    \vspace{.15in}
    \includegraphics[width=0.8\textwidth]{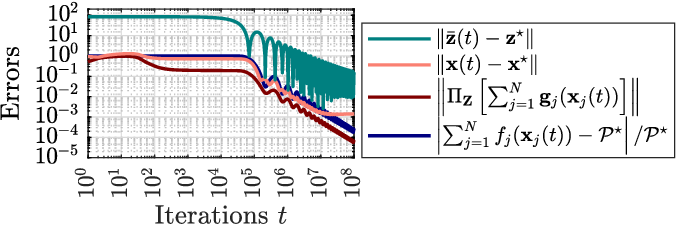}  
    \caption{Progress of Algorithm \ref{alg:ddsa} on $\Pcal_3$.
    }     
    \label{fig:td.2}	
\end{figure}

We report numerical results on a 204-bus T\&D system that comprises the IEEE 6-bus transmission network joined with six IEEE 33-bus distribution systems (see Figure \ref{fig:td.1} and \cite[Appendix B3]{extended} for details). We applied Algorithm \ref{alg:ddsa} on a reformulation of $\Pcal_3$ as an instance of $\Pcal$ with a flat start ($\v{z}_j(1) = 0$, $\v{Z}_j(0) = 0$, $j = 1, \ldots, N$) and step size $\eta = \eta_0/\sqrt{T}$, where $\eta_0 = 10^{3}$ and $T=10^{8}$. The agent-wise subproblems for $\Pcal_3$ are communicated over a 7-node star graph $\Gfk$ with the SO in the center. Convergence results are shown in Figure \ref{fig:td.2}. 

One might surmise that dual subgradient methods can be accelerated \`a la Nesterov. Such acceleration requires smoothness of the dual function--a property that an application problem such as $\Pcal_3$ does not always possess. In Appendix \ref{sec:dual_func}, we provide simple examples, where dual functions are nonsmooth, and show how the nonsmoothness impedes acceleration. We illustrate the same difficulty with $\Pcal_3$. 
Specifically, we compare the performance of Algorithm \ref{alg:ddsa} with a distributed acceleration scheme in \cite[eq. (5)-(6)]{qu2019accelerated}, adopted to our setup, described in Algorithm \ref{alg:Acc-ddsm_ges}. 
\begin{algorithm}
	\caption{Accelerated distributed dual subgradient to solve $\Pcal$.}
	\label{alg:Acc-ddsm_ges}

    Choose $\eta = \eta_0/\sqrt{T}$, 
    $\v{z}_j(1) = \v{Y}_j(1) = 0$, 
    $\alpha(1) = 0.5$, 
    $\v{x}_j(1) \gets \argmin_{\v{x}_j \in \Xset_j} \ \Lcal_j(\v{x}_j, \v{z}_j(1))$,
    and
    $\v{s}_j(1) \gets \v{g}_j(\v{x}_j(1))$.

	\For{$t=1, \ldots, T$}{	   
        
        $\v{Z}_j(t+1) \gets \sum_{k=1}^N W_{jk} \v{z}_k(t) + \eta \v{s}_j(t)$
        \label{alg:Acc-Z}
    
        $\v{Y}_j(t+1) \gets \sum_{k=1}^N W_{jk} \v{Y}_k(t)+ \frac{\eta}{\alpha(t)} \v{s}_j(t)$
        \label{alg:Acc-Y}

        $\alpha(t+1) \gets \frac{1}{2} \left[ -\alpha(t)^2 + \sqrt{\alpha(t)^4 + 4 \alpha(t)^2}  \right]$
        \label{alg:Acc-alpha}	

        $\v{z}_j({t+1})\gets \proj_\Zset \left[ ( 1 - \alpha(t+1) ) \v{Z}_j(t+1)\right. $ $\left. + \alpha(t+1) \v{Y}_j(t+1) \right]$
        \label{alg:Acc-z}	

	    $\v{x}_j(t+1) \gets \argmin_{\v{x}_j \in \Xset_j} \ \Lcal_j(\v{x}_j, \v{z}_j(t+1))$
        \label{alg:Acc-x}          

        $\v{s}_j({t+1})\gets  \sum_{k=1}^N W_{jk} \v{s}_k(t) 
            + \v{g}_j(\v{x}_j(t + 1)) -  \v{g}_j(\v{x}_j(t))$
        \label{alg:Acc-s}  
             
    }
\end{algorithm}

{
In this algorithm, the dual ancillary sequences $\v{Z}_j(t)$, $\v{Y}_j(t)$ in steps \ref{alg:Acc-Z}-\ref{alg:Acc-Y} are averaged by $\alpha(t)$ in step \ref{alg:Acc-z} that is updated in step \ref{alg:Acc-alpha}.
The distributed acceleration scheme in \cite[eq. (5)-(6)]{qu2019accelerated} has been designed for unconstrained optimization problems; we project multipliers of constraints \eqref{eq:TD.8} in step \ref{alg:Acc-z} on the positive orthant to apply it to $\Pcal_3$. 
Each local problem is solved in step \ref{alg:Acc-x} to obtain the local dual subgradient $\v{g}_j(\v{x}_j(t + 1))$. The central subgradient is approximated in step \ref{alg:Acc-s} that combines consensus tracking and local gradient averaging.

The results are shown in Figure \ref{fig:td.3}. Here, we adopted cold start for the dual iterates for both algorithms, i.e., $\v{z}_j(1) = \v{0},\ j = 1, \ldots, N$. We show the case when  generation costs for real power are considered quadratic on the left. The figure on the right is derived with linear real power generation costs. All reactive generation costs are considered quadratic for both simulations. 
Quadratic costs often yield a smooth dual function that Algorithm \ref{alg:Acc-ddsm_ges} can exploit and attain faster convergence than Algorithm \ref{alg:ddsa}. 
With linear costs, the dual function can be non-differentiable (especially at an optimal dual solution) as a result of primal non-uniqueness \cite[Theorem 6.3.3]{Bazaraa2006-oz}. 
Nonsmoothness is a fundamental barrier to acceleration as our example reveals. We  discuss the role of smoothness of the dual function in acceleration further in Appendix \ref{sec:dual_func} through illustrative examples.

\begin{figure}[htbp]
	\centering 
	\vspace{-0.05in}
	\includegraphics[width=0.85\textwidth]{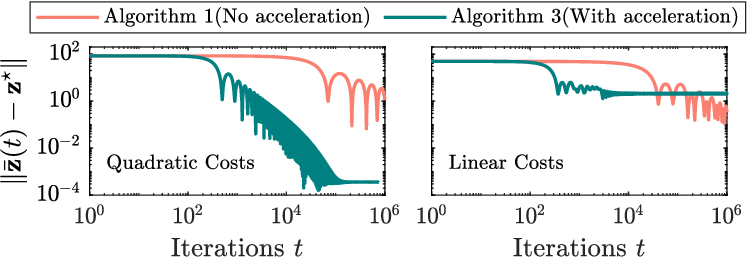}  
	\caption{{Convergence of dual iterates on $\Pcal_3$ about Algorithm \ref{alg:ddsa} (no acceleration) 
    and Algorithm \ref{alg:Acc-ddsm_ges} from \cite[eq. (5)-(6)]{qu2019accelerated} (with acceleration).}
	}     
	\label{fig:td.3}	
\end{figure}	
}

%% file: supplementary_materials.tex

{
\section{Dual subgradient methods cannot generally be accelerated}
\label{sec:dual_func}
Nesterov-type acceleration relies on smoothness of the objective function. Here, we illustrate why such acceleration is generally untenable in dual subgradient settings. We first present examples where the dual function for constrained optimization problems are nonsmooth. Then, we illustrate through examples how nonsmoothness of the objective function impairs acceleration.
Consider an optimization problem of the form
\begin{alignat}{2}
\underset{\v{x} \in \Xset}{\text{minimize}} \ 
    \frac{1}{2}\v{x}^\T\v{\Xi}\v{x} + \v{\xi}^\T\v{x} + \nu, 
    \ \ 
    \text{subject to} \ 
    \v{A} \v{x} \leq \v{b}.
    \label{eq:SP}
\end{alignat}
with $\v{\Xi}$ being positive semidefinite (but not positive definite) and $\Xset$ being a convex polyhedral set. This is a quadratic program (QP) that simplifies to a linear program (LP) when $\v{\Xi}=0$. Associate multipliers $\v{z} \geq 0$ with $\v{A}\v{x} \leq \v{b}$. Then, the dual function for the problem in \eqref{eq:SP} is given by
\begin{alignat}{2}
\begin{aligned}
    \Dcal(\v{z}):= 
    -\v{b}^\T \v{z} + \nu 
    + \underset{\v{x} \in \Xset}{\text{minimum}} \quad 
    \left\{ \frac{1}{2}\v{x}^\T\v{\Xi}\v{x} + \v{\xi}^\T \v{x} + \v{z}^\T \v{A} \v{x} \right\}.
    \end{aligned}
    \label{eq:dual_SP}
\end{alignat}
With parameter choices
\begin{gather}
\begin{gathered}
    \Xset := [0, 0.1]^{3},  
    \;
    \v{\xi} = -\begin{bmatrix}
        17 & 17 & 11 \\       
    \end{bmatrix}^\T, 
    \;
    \nu = 5, 
    \\
    \v{A} = \begin{bmatrix}
        0.19  &  0.12  &  0.42 \\
        0.37  &  0.54  &  0.13 \\
    \end{bmatrix}, 
    \;
    \v{b} = \begin{bmatrix}
        0.04 \\
        0.06 \\
    \end{bmatrix}
    \end{gathered}
    \label{eq:SP_coef}
\end{gather}
we plot $\Dcal(\v{z})$ with $\v{\Xi} = 0$ and $\v{\Xi} = \diag\left(24, 26, 0  \right)$, respectively, in the left and right of Figure \ref{fig:nonsmooth}. The primal and dual optimum is $\Pcal^\star = \Pcal^\star_{D} = 2.43$ for the QP and $2.30$ for the LP. 
\begin{figure}[htbp]
    \centering 
        \includegraphics[width=0.6\textwidth]{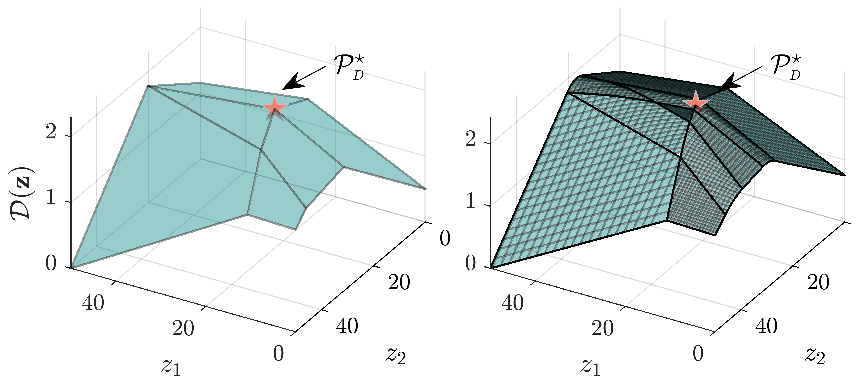}
    \\
    \includegraphics[width=0.6\textwidth]{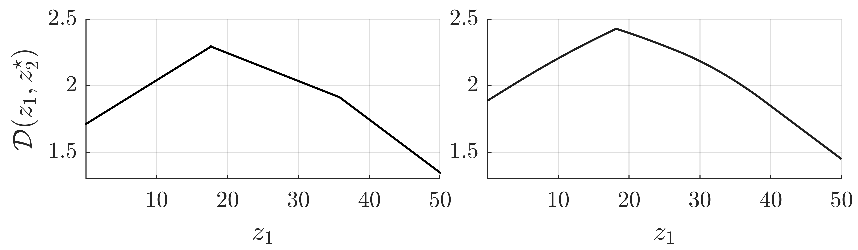}
    \caption{Illustrations of nonsmooth $\Dcal(\v{z})$ on top and a slice $\Dcal(\cdot, z_2^\star)$ on the bottom for \eqref{eq:SP} with parameters in \eqref{eq:SP_coef}. The plots on the left are obtained with $\v{\Xi} = 0$ and on the right with $\v{\Xi} = \diag\left(24, 26, 0  \right)$.}
    \label{fig:nonsmooth}
\end{figure}
$\Dcal$ in Figure \ref{fig:nonsmooth} is nonsmooth.
One might expect that for QP with $\v{\Xi} \neq 0$, $\Dcal$ might be smooth. Indeed with positive definite $\v{\Xi}$, $\Dcal$ becomes smooth, per \cite[Lemma 2.2]{necoaraLinearConvergenceDistributed2015a} or \cite[Theorem 6.3.3]{Bazaraa2006-oz}. Smoothness can no longer be guaranteed, when $\v{\Xi}$ is only positive semidefinite, as evidenced by Figure \ref{fig:nonsmooth}.


\newcommand{\xc}{{\v{x}_{\textrm{C}}}}
To demonstrate the impact of nonsmoothness on acceleration, consider the problem
    \begin{alignat}{2}
    \underset{\v{x} \in \Rset^{2}}{\text{minimize}} & 
    \quad f(\v{x}) := 
    \frac{1}{2} \| \v{x} - \xc \|_2^2
    + \lambda \| \v{x} \|_1.
    \label{eq:nonsmooth_compare}
    \end{alignat}
With $\xc = 0$, $f$ is minimized at the origin. If $\lambda =0$, then $f$ is smooth. For $\lambda>0$, $f$ is nonsmooth at the origin. With $\xc \neq 0$, $f$ is again smooth with $\lambda = 0$. However, with $\lambda > 0$, $f$ is nonsmooth at the origin, but not at the optimum of $f$. 



	
	



{
We compare (sub)gradient descent, i.e., $\v{x}(t+1) = \v{x}(t) - \eta \nabla f(\v{x}(t))$, with an accelerated variant described in Algorithm \ref{alg:Acc-Nesterov1983} on \eqref{eq:nonsmooth_compare}.
\begin{algorithm}[!ht]
	\caption{{Accelerated Gradient Descent adopted from \cite[Section 3.7.2]{bubeckConvexOptimizationAlgorithms2015}.}}
	\label{alg:Acc-Nesterov1983}

    Choose $\eta = 0.003$, 
    $\v{x}(1) = \v{y}(1) = (0.05, 0.05)$,
    and $\alpha(1) = 0$.

	\For{$t=1, \ldots, T$}{	   

        $\alpha(t+1) \gets \frac{1 + \sqrt{1 + 4 \alpha(t)^2}}{2}$, \  
        $\gamma(t) \gets \frac{1 - \alpha(t)}{\alpha(t+1)}$
        \label{alg:Acc-NGD.weight}           
        
        $\v{y}(t+1) \gets \v{x}(t) - \eta \nabla f(\v{x}(t))$
        \label{alg:Acc-NGD.y}

        $\v{x}(t+1) \gets ( 1 - \gamma(t) ) \v{y}(t+1) + \gamma(t) \v{y}(t)$
        \label{alg:Acc-NGD.x}

    }
\end{algorithm}

We start both algorithms at $(0.05, 0.05)$ with $\eta = 0.003$. For AGD, we also set $\alpha(1) = 0$.}
The left of Figure \ref{fig:nonsmooth_compare} shows the progress of subgradient descent (SGD) and accelerated (sub)gradient descent (AGD) on \eqref{eq:nonsmooth_compare} with $\xc = 0$. 
Note that when $\lambda = 0$, i.e., when $f$ is smooth, AGD outperforms SGD. However, when $\lambda = 0.01$ and $f$ is nonsmooth at the optimum, AGD performs better initially, but both algorithms oscillate around the optimum with similar errors, eventually.

In the right of Figure \ref{fig:nonsmooth_compare}, we compare SGD and AGD under the same settings, but with $\xc = (0.1,0.1)$. Note that AGD now performs better than SGD with zero and nonzero $\lambda$. That is, when the nonsmoothness is away from the optimum, AGD can accelerate convergence to the optimum locally, as long as the iterates remain within a region around the optimum where the function is locally smooth.
\begin{figure}[htbp]
    \centering 
    \includegraphics[width=0.6\textwidth]{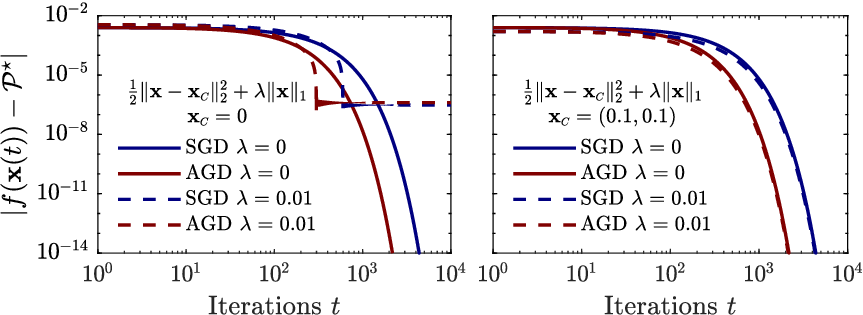}  
    \caption{{Comparison of SGD and AGD on \eqref{eq:nonsmooth_compare} with $\xc=0$ on the left and with $\xc=(0.1, 0.1)$ on the right.}}     
    \label{fig:nonsmooth_compare}
\end{figure}

}